\documentclass[12pt]{article}

\usepackage{amsmath, amssymb, eucal, amscd, color, amsthm, latexsym}

\begin{document}

\newtheorem{proposition}{Proposition}
\renewcommand{\theproposition}{}
\newtheorem{lemma}{Lemma}
\renewcommand{\thelemma}{}
\newtheorem{corollary}{Corollary}
\newtheorem{remark}{Remark}
\newtheorem{definition}{Definition}
\newtheorem{example}{Example}

\newcommand{\ca}{\mathcal}
\newcommand{\CC}{{\mathbb{C}}}
\newcommand{\PP}{{\mathbb{P}}}
\newcommand{\QQ}{{\mathbb{Q}}}
\newcommand{\RR}{{\mathbb{R}}}
\newcommand{\ZZ}{{\mathbb{Z}}}
\newcommand\Hom{\operatorname{Hom}}
\newcommand\End{\operatorname{End}}
\newcommand\Spec{\operatorname{Spec}}
\newcommand\Sym{\operatorname{Sym}}
\newcommand\CH{\operatorname{CH}}
\newcommand\ch{\operatorname{ch}}
\newcommand\Cone{\operatorname{Cone}}
\newcommand{\Tot}{\operatorname{Tot}}
\newcommand{\Ker}{\operatorname{Ker}}
\newcommand{\Cok}{\operatorname{Cok}}
\newcommand\Image{\operatorname{Im}}
\newcommand\sgn{\operatorname{sgn}}
\newcommand\Max{\operatornamewithlimits{Max}}

\newcommand{\cA}{{\cal A}} \newcommand{\cB}{{\cal B}}\newcommand{\cC}{{\cal C}}
\newcommand{\cD}{{\cal D}}\newcommand{\cE}{{\cal E}}\newcommand{\cF}{{\cal F}}
\newcommand{\cG}{{\cal G}}\newcommand{\cH}{{\cal H}}\newcommand{\cJ}{{\cal J}}
\newcommand{\cK}{{\ca K}} 
\newcommand{\cL}{{\ca L}} \newcommand{\cM}{{\ca M}} \newcommand{\cN}{{\ca N}}
\newcommand{\cR}{{\ca R}}
\newcommand{\cU}{{\cal U}} \newcommand{\cV}{{\cal V}}
\newcommand{\cW}{{\cal W}}  
 \newcommand{\cZ}{{\cal Z}}\newcommand{\Bf}{\mathbf{f} }

\newcommand{\bF}{{\mathbf{F}}}
\newcommand{\bU}{ {\mathbf{U}} }
\newcommand{\bK}{ {\mathbf{K}} }
\newcommand{\bL}{ {\mathbf{L}} }
\newcommand{\Symb}{\mathop{Symb}}
\newcommand{\bound}{\partial}

\newcommand{\chZ}{\check{\cal Z}}

\newcommand{\al}{\alpha} \newcommand{\be}{\beta} \newcommand{\ga}{\gamma} 
\newcommand{\de}{\delta} \newcommand{\eps}{\epsilon} \newcommand{\la}{\lambda}
\newcommand{\Ga}{\Gamma} \newcommand{\De}{\Delta}

\newcommand{\BA}{{\mathbb{A}}} \newcommand{\BB}{{\mathbb{B}}}
\newcommand{\BC}{{\mathbb{C}}} \newcommand{\BD}{{\mathbb{D}}}
\newcommand{\BE}{{\mathbb{E}}} \newcommand{\BF}{{\mathbb{F}}}
\newcommand{\BG}{{\mathbb{G}}} \newcommand{\BH}{{\mathbb{H}}} 
\newcommand{\BI}{{\mathbb{I}}} 
\newcommand{\BJ}{{\mathbb{J}}} \newcommand{\BK}{{\mathbb{K}}} 
\newcommand{\BL}{{\mathbb{L}}} \newcommand{\BM}{{\mathbb{M}}} 
\newcommand{\BN}{{\mathbb{N}}} \newcommand{\BO}{{\mathbb{O}}} 
\newcommand{\BP}{{\mathbb{P}}} \newcommand{\BQ}{{\mathbb{Q}}}
\newcommand{\BR}{{\mathbb{R}}} \newcommand{\BS}{{\mathbb{S}}}
\newcommand{\BT}{{\mathbb{T}}} \newcommand{\BU}{{\mathbb{U}}} 
\newcommand{\BV}{{\mathbb{V}}} 
\newcommand{\BW}{{\mathbb{W}}} 

\newcommand{\sq}{\square} 
\newcommand{\ul}{\underline}
\newcommand{\ulm}{ {\ul{m}\,} } \newcommand{\uln}{ {\ul{n}\,} } 
\newcommand{\ot}{\otimes}
\newcommand{\scirc}{\circ} 
\newcommand{\hts}{\hat{\otimes}}
\newcommand{\ts}{\otimes}
\newcommand{\vphi}{{\varphi}}  
\newcommand{\ddd}{\cdots}
\newcommand{\hlim}{\operatornamewithlimits{holim}}
\newcommand{\bd}{\partial}
\newcommand{\sAb}{s Ab}
\newcommand{\stm}{\times}
\newcommand{\bop}{\mathop{\textstyle\bigoplus}}
\newcommand{\Cyc}{Z} \newcommand{\C}{\Cyc}
\newcommand{\injto}{\hookrightarrow}
\newcommand{\hooklongrightarrow}{\lhook\joinrel\longrightarrow}
\newcommand{\surjto}{\twoheadrightarrow}
\newcommand{\isoto}{\overset{\sim}{\to}}
\newcommand{\J}{\cJ}
\newcommand{\ctop}{\overset{\,\,\circ}}
\newcommand{\minus}{\,\backslash\,}

\newcommand{\CHM}{{C\!H\!\cal M}}
\newcommand{\mapright}[1]{%
  \smash{\mathop{%
    \hbox to 1cm{\rightarrowfill}}\limits^{#1} } } 
\newcommand{\mapr}{\mapright}
\newcommand{\smapr}[1]{%
  \smash{\mathop{%
    \hbox to 0.5cm{\rightarrowfill}}\limits^{#1} } } 
\newcommand{\maprb}[1]{%
  \smash{\mathop{%
    \hbox to 1cm{\rightarrowfill}}\limits_{#1} } } 
\newcommand{\mapleft}[1]{%
  \smash{\mathop{%
    \hbox to 1cm{\leftarrowfill}}\limits^{#1} } }
\newcommand{\mapl}{\mapleft}
\newcommand{\maplb}[1]{%
  \smash{\mathop{%
    \hbox to 1cm{\leftarrowfill}}\limits_{#1} } }
\newcommand{\mapdown}[1]{\Big\downarrow
  \llap{$\vcenter{\hbox{$\scriptstyle#1\, $}}$ }}
\newcommand{\mapd}{\mapdown}
\newcommand{\mapdownr}[1]{\Big\downarrow
  \rlap{$\vcenter{\hbox{$\scriptstyle#1\, $}}$ }}
\newcommand{\mapdr}{\mapdownr}
\newcommand{\mapup}[1]{\Big\uparrow
  \llap{$\vcenter{\hbox{$\scriptstyle#1\, $}}$ }}
 \newcommand{\mapu}{\mapup}
\newcommand{\mapupr}[1]{\Big\uparrow
  \rlap{$\vcenter{\hbox{$\scriptstyle#1\, $}}$ }}
 \newcommand{\mapur}{\mapupr}

\newcommand{\Z}{\cZ} 
\newcommand{\ti}{\tilde} 
\newcommand{\La}{\Lambda} 
\newcommand{\Gr}{\operatorname{Gr}}
\def\dual{^{\vee} }

\title{ Integrals of logarithmic forms on semi-algebraic sets 
and a generalized Cauchy formula\\
\qquad Part I: convergence theorems}
\author{Masaki Hanamura, Kenichiro Kimura, and Tomohide Terasoma}
\date{}
\maketitle

\begin{abstract} 
In this paper consisting of two parts,
we study the integral of a logarithmic differential form on a compact semi-algebraic set in $\RR^n$ or $\CC^n$.
In Part I, we prove the convergence of the integral when the semi-algebraic set satisfies {\it  allowability} (or {\it admissibility}), a
 condition on the dimension of the intersection of the set and the pole
divisor of the differential form.  
\end{abstract}

\renewcommand{\thefootnote}{\fnsymbol{footnote}}

\newcommand{\Part}{Part\,}
\newcommand{\APropboundintegral}{2.5}
\newcommand{\sectbasic}{1}
\newcommand{\Propintegralposchart}{\sectbasic.1}
\newcommand{\Propintegralchart}{\sectbasic.2}
\newcommand{\PropIntegralinequality}{\sectbasic.3}
\newcommand{\Propmeasurezero}{\sectbasic.4}
\newcommand{\Propsumphi}{\sectbasic.5}
\newcommand{\Propaphi}{\sectbasic.6}
\newcommand{\Notationintegral}{\sectbasic.7}
\newcommand{\Nonoriented}{\sectbasic.8}
\newcommand{\Cvaluedforms}{\sectbasic.9}
\newcommand{\PropSreg}{\sectbasic.10}
\newcommand{\PropHomeonash}{\sectbasic.11}
\newcommand{\PropIsopartdense}{\sectbasic.12}
\newcommand{\CorPropIsopartdense}{\sectbasic.12.1}
\newcommand{\DefIntegralSalgset}{\sectbasic.13}
\newcommand{\Integralinvariance}{\sectbasic.14}

\newcommand{\sectResult}{2}
\newcommand{\Pullback}{\sectResult.1}
\newcommand{\Defdelta}{\sectResult.2}
\newcommand{\LemmaPartition}{\sectResult.3}
\newcommand{\PropPartitionforms}{\sectResult.4}
\newcommand{\Propboundintegral}{\sectResult.5} 
\newcommand{\ThmNashfinite}{\sectResult.6}
\newcommand{\PropBoundingAdd}{\sectResult.7}
\newcommand{\Propcontinuity}{\sectResult.8}
\newcommand{\Thmstokes}{\sectResult.9}

\newcommand{\sectAllow}{3}
\newcommand{\Defallowable}{\sectAllow.1}
\newcommand{\DefAllowable}{\sectAllow.2}
\newcommand{\basicproperties}{\sectAllow.3}
\newcommand{\lemmabasic}{\sectAllow.4}
\newcommand{\blochthm}{\sectAllow.5}
\newcommand{\propreduction}{\sectAllow.6}
\newcommand{\thmreduction}{\sectAllow.7}
\newcommand{\propfiniteness}{\sectAllow.8}
\newcommand{\PropVariation}{\sectAllow.9}
\newcommand{\PropLtypeineq}{\sectAllow.10}
\newcommand{\PropBlupfinite}{\sectAllow.11}
\newcommand{\Thmconvergence}{\sectAllow.12}
\newcommand{\Propnormallink}{\sectAllow.13}
\newcommand{\Thmnormallink}{\sectAllow.14}

\newcommand{\sectadmiss}{4}
\newcommand{\Dadmissible}{\sectadmiss.1}
\newcommand{\Spolar}{\sectadmiss.2}
\newcommand{\admissallow}{\sectadmiss.3}
\newcommand{\convadmiss}{\sectadmiss.4}
\newcommand{\admissallowext}{\sectadmiss.5}
\newcommand{\admissallowextD}{\sectadmiss.6}
\newcommand{\convadmisslink}{\sectadmiss.7}
\newcommand{\convadmisslinkD}{\sectadmiss.8}

{\bf Contents}\smallskip 

1. Integration on semi-algebraic sets

2. Basic results on integration of semi-algebraic sets

3. Integrals on allowable semi-algebraic sets

4. Integrals on admissible semi-algebraic sets

\bigskip

This paper consists of two parts.
The purpose of \Part I  is the study of integrals of differential forms 
on semi-algebraic sets on $\RR^n$ or $\CC^n$.
In particular, we investigate conditions under which the integral converges.
In \Part II, we formulate and prove a generalization of the Cauchy formula
for integrals on semi-algebraic sets. 
The formula will play a key role in constructing the ``Hodge realization" functor from the
category of mixed Tate motives over $\CC$. 

The integral we consider is typically of the form
$$\int_A \frac{dz_1}{z_1}\wedge\cdots \wedge \frac{dz_n}{z_n}$$
where $A$ is a compact semi-algebraic set 
 of dimension $n$ in $\CC^n$.
For instance,  for a real number $a$ with $0<a<1$,  we have
$\int_{[a, 1]} dz/z= -\log (a)$; 
 for the set $S_a=\{(t_1, t_2)\in \RR^2\,|\, 0\le t_1\le 1, 0\le t_2\le a, t_1+t_2\ge 1\}$,  we have
$$\int_{S_a}\frac{dz_1}{z_1}\wedge \frac{dz_2}{z_2} = -Li_2 (a)\,$$
where $Li_2(z)=\sum_{n\ge 1} z^n/n^2$, $|z|<1$, is the dilogarithm function.
More generally values of the polylogarithm functions $Li_p(z)=\sum_{n\ge 1} z^n/n^p$, $p\ge 2$, 
can be expressed as integrals on $\CC^n$. 
Other examples arise as the period integrals of the cohomology $H^n(\PP^n-(\cup \BH_i), X- (\cup \BH_i)\,)$,
where $\BH_i=\{z_i=0\}$, $i=0, \cdots, n$,  are the coordinate hyperplanes and $X$ is a subvariety of $\PP^n$. 

For examples of such integrals which naturally appear in algebraic geometry, the convergence may be verified by direct methods.
On the other hand, one may seek for a general geometric condition for the convergence.
Answering this question, the results in this paper, Theorem (\Thmconvergence) and  Theorem (\convadmiss), 
ensure the convergence under a simple assumption on the dimension of the intersection of the semi-algebraic set with the pole divisors of the form.
\bigskip

In \S 1, we give a general discussion on  the definition and properties of integrals
of differential forms on a semi-algebraic set, 
based on the Lebesgue integration theory.
For this we first describe Lebesgue integration of a smooth form on a smooth manifold;
this material is included for lack of appropriate reference and for setting up the necessary 
notation.
We then study integration on a semi-algebraic set.
Let $S$ be an oriented closed semi-algebraic set of $\RR^n$, with $\dim S=m$, 
and let $\vphi$ be a smooth $m$-form defined on an open semi-algebraic set
$U$ of $S$, which is a Nash submanifold of $\RR^n$, with $\dim (S\backslash U)<m$. 
Then one can define the integral $\int_S|\vphi|$ as $\int_U|\vphi|$, which is a non-negative
real number (possibly infinite). 
When $\int_S|\vphi|<+\infty$ (we then say that $\int_S \vphi$ is absolutely convergent)
 one can also define the integral $\int_S \vphi$, which is a real number.
For the arguments in \S 1, some basic results  on semi-algebraic sets are needed;  we have indicated how
they can be derived from known facts
found in [BCR].

In \S 2, we establish some basic results which will be useful
throughout our paper.
We first show Proposition (\APropboundintegral).  
As a consequence, if $S$ is compact semi-algebraic of dimension $m$, 
$h: S\to \RR^\ell$ a continuous semi-algebraic map, and $\psi$ a smooth $m$-form 
defined on a neighborhood of $h(S)$ in $\RR^\ell$, then 
the integral $\int_S |h^*\psi|$ is finite. (This is an ``easy" convergence theorem.)

Another consequence of Proposition (\APropboundintegral) is the Stokes formula:
 if $S\subset \RR^n$ is a semi-algebraic 
set, $S=\Delta^m$ as a semi-algebraic set of $\RR^n$, and $\psi$ a smooth $(m-1)$-form defined
on a neighborhood of $S$, then one has
$$\int_{\delta S}  \psi =\int_{ S} d\psi\,.$$
Note that it has been known that the Stokes formula 
holds more generally for sub-analytic sets (see [He]).
See also [OS] for a different approach in the case of semi-algebraic sets.

In \S 3, which is central to \Part I,
 we consider the following question:
Given an oriented compact semi-algebraic set $A\subset \RR^n$ of dimension $n$
and a smooth $n$-form $\omega$ defined on an open set of $A$, under 
what condition the integral $\int_A \omega$ converges?

We assume that $\omega$ has only logarithmic singularities along the 
coordinate hyperplanes.  An non-empty intersection of coordinate hyperplanes
is called a face.
We introduce the notion of {\it allowability} for a closed semi-algebraic set $A$ of $\RR^n$ as follows:  $A$ is allowable if the intersection
of $A$ with any face $F$ has dimension strictly less than $\dim F$. 
Our main result is Theorem (\Thmconvergence):  If $A$ is compact and allowable in 
$\RR^n$, and $\omega$ is a logarithmic $n$-form, then the integral 
$\int_A \omega$ absolutely converges. 
The proof is outlined at the beginning of \S 3.

In addition, we consider the slice of $A$ by a coordinate of $\RR^n$, and study its 
volume with respect to a logarithmic form.
More precisely, for $A$ compact and allowable, let  
$A_t=A\cap \{x_1=t\}$, the intersection of $A$ with
the hyperplane $\{x_1=t\}$ for $t\in \RR$ with $|t|$ small and non-zero.
Since $A_t$ is of dimension $\le n-1$, if $\omega$ is an $(n-1)$-form on $\RR^n$
with only logarithmic singularities, one has the integral  $\int_{A_t}|\omega|$; 
we show in Theorem (\Thmnormallink) that it tends to zero as $t\to 0$.

In \S 4, we study integrals on semi-algebraic sets of $\CC^n$. 
For a compact semi-algebraic set $A$ of $\CC^n$ of dimension
$m$ (where $n\le m\le 2n$), and for an $(n, m-n)$-form $\omega$ on $\CC^n$ with only
logarithmic singularities along the coordinate hyperplanes,
we ask when the integral $\int_A \omega$ absolutely converges.
We define the notion of {\it $m$-admissibility} on $A$, and 
when $A$ is $m$-admissible, we show that the integral converges, Theorem (\convadmiss). 
The proof consists of reduction to an integral of an allowable algebraic set in $\RR^m$
by means of a change of coordinates (essentially by the polar coordinates) and then applying 
Theorem (\Thmconvergence). 
Also, as a counterpart of Theorem (\Thmnormallink), 
we show in Theorem (\convadmisslink) that, for $A$ compact and $m$-admissible, 
if $A_t :=A\cap \{|z_1|=t \ge |z_2|\}$ for $t\in\RR$ with $|t|$ small, and 
$\omega$ is an $(n, m-n-1)$-form on $\CC^n$ with 
logarithmic singularities, then 
$\int_{A_t} |\omega|$ tends to zero as $t\to 0$. 
Theorems (\convadmiss) and (\convadmisslink) will be effectively used in the proof of 
of the Cauchy formula in \Part II.
\bigskip 

We would like to thank Professor T. Suwa for helpful discussions on
the intersection theory of semi-algebraic sets; they are incorporated in \Part II.
\bigskip

\newpage

\newcommand{\RRinfty}{\RR_{\ge 0}\cup\{+\infty\}}

{\bf \S \sectbasic. Integration on semi-algebraic sets}\bigskip 

{\bf Integration on a manifold.}\quad
Let $M$ be an $m$-dimensional oriented $C^\infty$-manifold
(satisfying the second axiom of countability).
By an $m$-form $\vphi$ on $M$, we mean a possibly discontinuous
section of $\wedge^m T^*M$, the $m$-th wedge product of the cotangent bundle of $M$;
thus on each oriented chart $U$ with coordinates $(x_1, \cdots, x_m)$, one can write 
$$\vphi|_U=g\,dx_1\wedge\cdots\wedge dx_m$$
where $g$ is a function on $U$ with values in $\RR$. 
We say that $\vphi$ is a measurable $m$-form if, on each chart $U$,  $g$ is 
Borel measurable function. 

For convenience of the reader, we recall that the $\sigma$-algebra $\cB$ of 
the Borel measurable sets of $M$ is the smallest $\sigma$-algebra containing all open sets.
A function $f$ defined on $E\in \cB$ and takes values in $\RR\cup\{\pm \infty\}$ 
is Borel measurable if for any $c\in \RR\cup\{\pm \infty\}$, the set
$\{x\in E\, |\, f(x)>c\}$ is in $\cB$. 

A measurable $m$-form $\vphi$ is {\it non-negative} with respect to the orientation if on each oriented chart $U$, the function $g$ is non-negative.
Let $\{U_\al\}$ be a countable covering by oriented charts, and $\{\rho_\al\}$ be a partition of unity subordinate to the covering.
For each $U_\al$,  write
$\vphi|_{U_\al}= g_\al\, dx_1\wedge\cdots\wedge dx_m$ 
where $(x_1, \cdots, x_m)$ are the local coordinates and $g_\al=g_\al(x_1, 
\cdots, x_m)$ is a non-negative measurable function.
For a Borel measurable set $E$ of $M$, define $\mu_\vphi(E)\in \RR_{\ge 0} \cup \{+\infty\}$
as the countable sum
$$\mu_\vphi(E)=\sum_\al \int_{E\cap U_\al} \rho_\al  g_\al \,dx_1\cdots dx_m
 \,,$$
where each term $\int_{E\cap U_\al} \rho_\al  g_\al \,dx_1\cdots dx_m$
is the Lebesgue integral of the measurable function $\rho_\al g_\al$
on the Borel set $E\cap U_\al $ of $U_\al$; note that $U_\al$ is identified 
via the local coordinates with an open subset of $\RR^m$, where the usual Lebesgue
integration theory applies.
The value $\mu_\vphi(E)$ is well-defined, independent of the choice of a covering $\{U_\al\}$ and 
$\{\rho_\al\}$;  this can be shown  
using the change of variables formula for 
the Lebesgue integral.

The set function $\mu_\vphi(E)$ satisfies the complete additivity:
If $(E_k)_{k=1, 2, \cdots}$ is a disjoint family of Borel measurable sets, 
then 
$$\mu_\vphi(\sum_{k=1}^\infty E_k)=\sum_{k=1}^\infty\mu_\vphi( E_k)\,.$$
Indeed, setting $E= \sum_k E_k$, 
by the complete additivity of the Lebesgue integral
on each $U_\al$, one has
$$\int_{E\cap U_\al} \rho_\al g_\al \, dx_1\cdots dx_m
=\sum_k \int_{E_k\cap U_\al} \rho_\al g_\al \, dx_1\cdots dx_m\,.$$
Taking the sum over $\al$ gives the assertion.

We thus obtain a measure space $(M, \cB, \mu_\vphi)$, 
where $\cB$ is the $\sigma$-algebra of Borel measurable sets
of $M$.
For this measure space we have the usual 
Lebesgue integration theory, some of which we briefly recall:
\smallskip 

$\bullet$\quad For a non-negative Borel measurable function $f$ defined on a measurable set $E$, we can define the 
integral 
$$\int_E f \,d\mu_\vphi\in \RRinfty\,.$$
(The definition will be briefly recalled in the proof of Proposition (\Propintegralposchart). )
We say that $f$ is integrable with respect to $\mu_\vphi$ if $\int_M f \,d\mu_\vphi<+\infty$.

In what follows, we will consider only Borel measurable sets and Borel measurable functions, so we will often drop 
the prefix ``Borel".

$\bullet$\quad For $f$ a Borel measurable
function on $E$, one has $f=f^+ - f^-$, where $f^+(x)=Max (f(x), 0)$
and $f^-(x)=Max (-f(x), 0)$ are non-negative Borel measurable functions.
We say that $f$ is absolutely integrable with respect to $\mu_\vphi$ if $f^+$ and $f^-$ are both integrable
(equivalently if $|f|$ is integrable);
then we can define the integral of $f$ on $E$ by
$$\int_E f\, d\mu_\vphi =\int_E f^+ \,d\mu_\vphi  -\int_E f^- \,d\mu_\vphi\in \RR\,.$$
\bigskip

It is useful to note that the integral $\int f\, d\mu_\vphi$ can be expressed as a sum of integrals on charts:
\bigskip 

\noindent (\Propintegralposchart) {\bf Proposition.}\quad{\it  Let $M$, 
$\{U_\al\}$, and $\{\rho_\al\}$ 
be as above.  
Let $\vphi$ be a non-negative measurable $m$-form on $M$
and, and write $\vphi|_{U_\al}= g_\al \, dx_1\wedge\cdots\wedge dx_m$ with $g_\al$
a non-negative measurable function on each $U_\al$.

(1) For a non-negative Borel measurable function $f$
on a Borel measurable set $E$ of $M$, we have
$$\int_E f \, d\mu_\vphi =\sum_\al \int_{E\cap U_\al} \rho_\al  fg_\al \, dx_1\cdots dx_m\,\in \RRinfty\,.$$

(2) Let $f$ be a Borel measurable function on $E$, which is absolutely 
integrable with respect to $\mu_\vphi$.
Then, for each $\al$, the integral 
$\int_{E\cap U_\al} \rho_\al  fg_\al \, dx_1\cdots dx_m$ is absolutely convergent; further the sum 
$\sum_\al\int_{E\cap U_\al} \rho_\al  fg_\al \, dx_1\cdots dx_m$ is absolutely convergent and we have 
$$\int_E f \, d\mu_\vphi =\sum_\al \int_{E\cap U_\al} \rho_\al  fg_\al \, dx_1\cdots dx_m\,\in \RR\,.$$
}\smallskip 

\begin{proof}
(1) If $f=\chi_E$, the characteristic function of  a Borel measurable set $E$, 
then the identity amounts to
$$\mu_\vphi(E)=\sum_\al \int_{E\cap U_\al} \rho_\al 
g_\al \, dx_1\cdots dx_m
 \,,$$
which is the definition of $\mu_\vphi(E)$. 
If $f$ is a simple function, namely if $f=\sum c_i\chi_{E_i}$ for a disjoint finite family of Borel measurable sets $\{E_i\}$, the identity holds by linearity. 

For an arbitrary measurable function $f\ge 0$ on $E$, there exists a monotonically increasing sequence 
of simple functions $f_n(x)$ such that $\lim f_n(x) =f(x)$. 
For each $n$ we have the identity 
$$\int_E f_n \, d\mu_\vphi =\sum_\al \int_{E\cap U_\al} \rho_\al f_n g_\al \, dx_1\cdots dx_m\,.$$
When $n\to \infty$, the left hand side tends to $\int_E f \, d\mu_\vphi$
by the definition of the Lebesgue integral. 

On the right hand side, for each $\al$, $\rho_\al f_n g_\al$ is monotonically 
increasing and $\lim \rho_\al f_n g_\al = \rho_\al f g_\al$;
by a well-known property of the Lebesgue integral, we have 
$$\lim_{n\to \infty}\int_{E\cap U_\al} \rho_\al f_n g_\al \, dx_1\cdots dx_m
=\lim_{n\to \infty}\int_{E\cap U_\al} \rho_\al f g_\al \, dx_1\cdots dx_m\,.$$
The identity hence follows. 

(2) If $f$ is absolutely 
integrable, one can apply (1) to $f^+$ and $f^-$, so that 
$$\int_E f^{\pm} \, d\mu_\vphi =\sum_\al \int_{E\cap U_\al} \rho_\al  f^\pm g_\al \, dx_1\cdots dx_m\,\in \RR\,.$$
The assertion hence follows. 
\end{proof}

Next assume we are given a measurable $m$-form $\vphi$ on $M$. 
On each oriented chart $U$, one has 
$\vphi|_U=g\, dx_1\wedge\cdots\wedge dx_m$, 
where $g$ is a measurable function
on $U$. We have non-negative measurable $m$-forms 
$\vphi^+$ and $\vphi^-$ on $M$ with the property that  
$\vphi=\vphi^+ -\vphi^-$ and on an oriented chart, 
$$\vphi^+|_U=g^+dx_1\wedge\cdots\wedge dx_m \quad
\text{and}\quad \vphi^-|_U=g^-dx_1\wedge\cdots\wedge dx_m
\,$$
where $g=g^+-g^-$ is the decomposition of $g$ to the positive and negative parts. 
Then $|\vphi|:=\vphi^+ +\vphi^-$ is another non-negative measurable 
$m$-form, and on an oriented chart, 
$|\vphi|\, |_U=|g|dx_1\wedge\cdots\wedge dx_m$. 

We can apply the discussion before to the non-negative forms $\vphi^+$, $\vphi^-$ and 
$|\vphi|$; associated with these $m$-forms are
 the measures $\mu_{\vphi^+}$, $\mu_{\vphi^-}$
 and $\mu_{|\vphi|}=\mu_{\vphi^+} +\mu_{\vphi-}$.

Let now $f$ be a  Borel measurable function  on a Borel measurable set $E$. 
If $f$ is non-negative, we have the integral with respect to the measure $\mu_{|\vphi|}$, 
$$\int_E fd \mu_{|\vphi| }=\int_E fd \mu_{\vphi^+} +\int_E fd \mu_{\vphi^-}
\in \RRinfty.
$$
When $f$ is not necessarily non-negative, but absolutely integrable 
with respect to $\mu_{|\vphi| }$ (thus also with respect to $\mu_{\vphi^+}$ and $\mu_{\vphi^-}$), 
we have  
$$\int_E fd \mu_{\vphi^+}\in \RR\quad\text{and} \quad \int_E fd \mu_{\vphi^-} \in \RR\,;$$
thus we shall set 
$$\int_E fd \mu_{\vphi} =\int_E fd \mu_{\vphi^+} -\int_E fd \mu_{\vphi^-} \in \RR.$$

From  Proposition (\Propintegralposchart) and the definitions we obtain:
\bigskip 

\noindent (\Propintegralchart) {\bf Proposition.}\quad{\it  
Let $\vphi$ be a measurable $m$-form on $M$, and write
$\vphi|_{U_\al}= g_\al dx_1\wedge\cdots\wedge dx_m$ with a measurable function
$g_\al$
on each $U_\al$. 

(1) For a non-negative Borel measurable function $f$
on a Borel measurable set $E$ of $M$, we have
$$\int_E f \, d\mu_{\vphi^+} =\sum_\al \int_{E\cap U_\al} \rho_\al f  g^+_\al \,dx_1\cdots dx_m\in \RRinfty\,.$$
Similarly for $\vphi^-$ and $|\vphi|$, with $g_\al^+$ replaced by $g_\al^-$ and $|g_\al|$, respectively. 

(2) Let $f$ be a Borel measurable function on $E$, which is absolutely 
integrable with respect to $\mu_{|\vphi|}$.
Then, for each $\al$, the integral 
$\int_{E\cap U_\al} \rho_\al  fg_\al \, dx_1\cdots dx_m$ is 
absolutely convergent; further the sum 
$\sum_\al\int_{E\cap U_\al} \rho_\al  fg_\al \, dx_1\cdots dx_m$ is absolutely convergent and we have 
$$\int_E f \, d\mu_\vphi =\sum_\al \int_{E\cap U_\al} \rho_\al  fg_\al \, dx_1\cdots dx_m\,\in \RR\,.$$
}

\begin{proof}
(1) One applies (\Propintegralposchart), (1) to the form $\vphi^+$. 

(2) Writing $f=f^+-f^-$, apply (1) to $f^+$ and $f^-$ to obtain  four equalities 
$$\int_E f' \, d\mu_{\vphi''} =\sum_\al \int_{E\cap U_\al} \rho_\al f'g''_\al \, dx_1\cdots dx_m$$
where $'$ and $''$ are elements of $\{\pm\}$. 
By the assumption of absolute convergence, all the four values are finite, and we obtain the assertions.
\end{proof}

The following is obvious from the definitions.
\bigskip 

\noindent (\PropIntegralinequality) {\bf Proposition.}\quad{\it  
Let $\vphi$ be a measurable $m$-form on $M$, 
$E$ a Borel measurable set, and
$f$ be a Borel measurable function on $E$.
Then $f$ is absolutely integrable with respect to $\mu_{|\vphi|}$
if and only if $\int_E |f|\, d\mu_{|\vphi|}<+\infty$; we then have
$$\left|\int_Ef\, d\mu_{\vphi}\right|\le \int_E |f|\, d\mu_{|\vphi|}\,.$$
}\bigskip 

\noindent (\Propmeasurezero) {\bf Proposition.}\quad{\it  
A submanifold $N$ of $M$ with $\dim N <m$ is a set of measure zero with respect to $\mu_{|\vphi|}$. 
In particular, for any non-negative measurable function $f$ on $N$, the integral $\int_N f\,d\mu_{|\vphi|}=0$. 
}

\begin{proof}
The intersection $N\cap U_\al$ is a submanifold of $U_\al$ of  dimension $<m$, and it is 
easy to show that it has
measure zero with respect to the Lebesgue measure $dx_1\cdots dx_m$. 
Thus $\int_{E\cap U_\al} \rho_\al  |g_\al| \,dx_1\cdots dx_m=0$, and it follows from the definition of $\mu_{|\vphi|}$ that 
$$\mu_{|\vphi|}(N)=\sum_\al \int_{N\cap U_\al} \rho_\al |g_\al| \,dx_1\cdots dx_m=0\,.$$
\end{proof}

If $\vphi$ and $\psi$ are measurable $m$-forms on $M$, $\vphi +\psi$ is  a measurable 
$m$-form.   
\bigskip

\noindent (\Propsumphi) {\bf Proposition.}\quad{\it  
(1) Let $\vphi$ and $\psi$ be non-negative measurable $m$-forms on $M$. 
For a non-negative measurable function $f$ on $E$, we have 
$$\int_E f\, d\mu_{\vphi+\psi} =\int_E f\, d\mu_{\vphi}+\int_E f\, d\mu_{\psi}\in \RRinfty\,.$$

(2) Let $\vphi$ and $\psi$ be measurable $m$-forms on $M$. 
For a non-negative measurable function $f$ on $E$, we have an inequality 
$$\int_E f\, d\mu_{|\vphi+\psi|} \le \int_E f\, d\mu_{|\vphi|}+\int_E f\, d\mu_{|\psi|}\,.$$
If $f$ is a measurable function, which is absolutely integrable with respect to $\mu_{|\vphi|}$ and $\mu_{|\psi|}$, then 
$f$ is also absolutely integrable with respect to $\mu_{|\vphi+\psi|}$.

(3) When $f$ is absolutely integrable with respect to $\mu_{|\vphi|}$ and $\mu_{|\psi|}$, we have
$$\int_E f\, d\mu_{\vphi+\psi} =\int_E f\, d\mu_{\vphi}+\int_E f\, d\mu_{\psi} \,.$$
}\smallskip

\begin{proof}
(1) and (2) follow from (\Propintegralposchart), (1).

(3) The first assertion follows from (2) applied to $|f|$, and the second follows from (\Propintegralchart), (2).
\end{proof}
\bigskip

If $\vphi$ is an $m$-form on $M$ and $f$ a function defined on a subset $E$, 
then $f$ may be viewed as a function on $M$ (which takes the value zero outside
$E$), so $f\vphi$ is an $m$-form on $M$ which vanishes outside $E$. 
\bigskip 

\noindent (\Propaphi) {\bf Proposition.}\quad{\it  
(1) Let $\vphi$ be a non-negative measurable $m$-form on $M$, and $f$ be a non-negative $\RR$-valued measurable function on a measurable set $E$.
Then we have
$$\int_E f\,d\mu_{\vphi}= \int_E d\mu_{f\vphi }\in \RRinfty\,.$$

(2) Let $\vphi$ be a measurable $m$-form on $M$,  and $f$ be an $\RR$-valued measurable function on a  measurable set $E$. 
Then $f$ is absolutely convergent with respect to $\mu_{|\vphi|}$ if and only if 
$1$ is absolutely convergent with respect to $\mu_{|f\vphi|}$; in this case we have 
$$\int_E f\,  d\mu_{\vphi}=\int_E d\mu_{f\vphi} \in \RR\,.$$
}\smallskip 

\begin{proof}
(1)  Follows from (\Propintegralposchart), (1).

(2) The first assertion follows from (1) applied to $|f|$
and $|\vphi|$.
The second assertion follows from (\Propintegralchart), (2).
\end{proof}

This shows that little is lost if we consider only the case $f=1$.
Therefore, in setting up the following notation, we assume $f=1$.
\bigskip

\noindent (\Notationintegral) {\bf Notation.}\quad 
Let $M$ be an oriented manifold of dimension $m$. 
For a measurable $m$-form $\vphi$ on $M$ 
on a measurable set $E$, we will write 
$$\int_E \vphi\,\quad \text{for}\quad \int_E d\mu_\vphi\,, $$
and similarly 
$$\int_E  |\vphi| \quad \text{for}\quad \int_E d\mu_{|\vphi|}\,.$$

With this notation, we have
$$\int_E fd\mu_\vphi = \int_E f \vphi\,$$
by (\Propaphi), (2); also we have
$$\int_E  |\vphi+\psi| \le \int_E  |\vphi| +\int_E  |\psi|\,$$
by (\Propsumphi), (2), and 
$$\int_E (\vphi+\psi)= \int_E \vphi+ \int_E \psi\,$$
by (\Propsumphi), (3).
\bigskip 

\noindent (\Nonoriented) {\it Unoriented case.}\quad
Assume $M$ is not oriented (or even not orientable).  Let 
$\{U_\al\}$ be a countable covering by unoriented charts, and $\{\rho_\al\}$ be a partition
 of unity subordinate to the covering. 
Given a measurable $m$-form $\vphi$ on $M$, on each 
 $U_\al$, write $\vphi|_{U_\al}= g_\al dx_1\wedge\cdots\wedge dx_m$ 
where $(x_1, \cdots, x_m)$ are the local coordinates and $g_\al=g_\al(x_1, 
\cdots, x_m)$ is a measurable function.
 For a Borel measurable set $E$ of $M$, we set 
$$\mu_{|\vphi|}(E)=\sum_\al \int_{E\cap U_\al} \rho_\al  |g_\al| dx_1\cdots dx_m\in \RRinfty
 \,;$$
 then $\mu_{|\vphi|}$ defines a measure on $(M, \cB)$. 
In case $M$ is oriented, this coincides with the measure $\mu_{|\vphi|}$ defined earlier, which 
justifies the use of the same notation.  In general, however, note that $|\vphi|$  does not 
have a meaning.

If $f$ is a non-negative Borel measurable function taking values in $\RRinfty$  defined on a measurable set $E$, we 
have the 
integral 
$$\int_E f \,d\mu_{|\vphi|}\in \RRinfty\,.$$
If $f$ is a measurable function on $E$, which is absolutely integrable, the integral 
$\int_E f \,d\mu_{|\vphi|}\in \RR$ is defined. 

We have an analogue of (\Propintegralposchart), where $g_\al$ should be replaced with $|g_\al|$; the proof is the same. 
We also have analogues of  (\Propmeasurezero),  (2) of (\Propsumphi),  and (\Propaphi).

As in the oriented case, we will  write $\int_E  |\vphi|$ for $\int_E d\mu_{|\vphi|}$.
\bigskip 

\noindent (\Cvaluedforms)\quad  By standard arguments one may generalize all of the definitions and
results to the case of $\CC$-valued forms $\vphi$ and $\CC$-valued functions
$f$.   For example, (\PropIntegralinequality), (\Propmeasurezero), 
(2) and (3) of (\Propsumphi),  and (2) of (\Propaphi) hold. 
\bigskip

{\bf Integration on a semi-algebraic set.}\quad 
Let $V$ be an algebraic set of $\RR^n$. 
We say that $x\in V$ is a non-singular point in dimension $d$ if there is 
an irreducible component $V'$ of $V$ with $\dim V'=d$, such that
$V'$ is the only irreducible component of $V$
 containing $x$, and $x$ is a
non-singular point of $V'$. (See [BCR, Definition 3.3.9 and Proposition 3.3.10].)

If $V$ is an algebraic set of dimension $d$, then the set of non-singular
points of $V$ in dimension $d$, denoted $Reg(V)$, is a Zariski open set of $V$
and the complement $V-Reg(V)$ is an algebraic subset of dimension smaller than
$d$, [BCR, Proposition 3.3.14].
Also, $Reg(V)$ is a Nash submanifold of dimension $d$ of $\RR^n$
([BCR, Proposition 3.3.11]); for the definition of a Nash submanifold
see [BCR,  \S 2.9].

Let $S\subset \RR^n$ be a semi-algebraic set of dimension $m$. 
Let $V$ be the Zariski closure of $S$ (so $V$ is of dimension $m$), 
and $\ctop{S}$ be the interior of $S$ in $V$.
\bigskip 

\noindent (\PropSreg) {\bf Proposition.}\quad{\it  The set 
$S_{reg}:=\ctop{S}\cap Reg(V)$ is a non-empty open semi-algebraic subset of 
$S$, which is a Nash submanifold of $\RR^n$ of dimension $m$.
}\smallskip 

\begin{proof}
The set $S':=S\cap Reg(V)$ is a non-empty semi-algebraic subset of dimension $m$,  since $\dim (V-Reg (V))<m$.
As $\ctop{S}$ is an open semi-algebraic set of $Reg(V)$, 
$S_{reg}$ is either empty or an open Nash submanifold of $V$.
It thus remains to show that $S_{reg}$ is not empty.

By [BCR, Proposition 2.9.10], $S'$ is a disjoint union of a finite number of Nash submanifolds $S'_i$ with dimension $d_i\le m$.
Take an $S'_i$ with $d_i=m$ (there exists at least one);
then it is an open Nash submanifold of $Reg(V)$, so it is contained in 
$S_{reg}$, and the assertion is proven.
\end{proof}

We call $S_{reg}$ 
{\it  the regular part} of $S$. 
Note that $\dim (S\, \backslash\, S_{reg})<m$. 
\bigskip

\noindent (\PropHomeonash)
{\bf Proposition.}\quad{\it Let $f: M\to N$ be a surjective Nash map between Nash submanifolds. 

(1) There is an open Nash submanifold $U$ of $N$ with $\dim (N\minus U) <\dim N$ such that 
$f|_{f^{-1}(U)}: f^{-1}(U)\to U$ is submersive (i.e. for each $x\in f^{-1}(U)$, 
the rank of the differential map $d_x f: T_x M\to T_{f(x)}N$ is equal to
$\dim N$). 

(2) Assume further that $\dim M=\dim N=m$.
Then the map $f: f^{-1}(U)\to U$ of (1) is a local Nash diffeomorphism
(i.e., for any point $x\in f^{-1}(U)$, there are open semi-algebraic neighborhoods $W$ of $x$ and $U'$ of $f(x)$ such 
that $f|_W$ is a Nash diffeomorphism from $W$ to $U'$).
In particular, if $f$ is a semi-algebraic homeomorphism, then $f:f^{-1}(U)\to U$ is a Nash diffeomorphism.
}\smallskip 

\begin{proof}
 (1) This follows from Sard's theorem for Nash maps [BCR, Theorem 9.6.2], 
which states that the set of critical values of $f$ is a semi-algebraic set of dimension smaller
than $\dim N$.  

(2)  The differential map is an isomorphism, so the claim follows from 
the semi-algebraic inverse function theorem [BCR, Proposition 2.9.7].
\end{proof}

\noindent (\PropIsopartdense)
{\bf Proposition.}\quad{\it  Let $f: S\to S'$ be a semi-algebraic homeomorphism of 
semi-algebraic sets of dimension $m$. 
Then there exist open semi-algebraic sets $U\subset S$, $U'\subset S'$,
with $\dim (S\backslash U)<m$, $\dim (S'\backslash U')<m$, such that 
$U$, $U'$ are Nash submanifolds and $f$ induces a Nash diffeomorphism
$U\to U'$.}
\smallskip 

\begin{proof} By (\PropSreg) one may assume that $S$ and $S'$ are Nash submanifolds.
We conclude by applying (\PropHomeonash), (2).
\end{proof}

\noindent (\CorPropIsopartdense)
{\bf Corollary.}\quad{\it 
Let $S$ be a semi-algebraic set of dimension $m$ and $f: S\to \RR$ be a continuous semi-algebraic function. 
Then there exists an open semi-algebraic set $U$ of $S$ with $\dim (S\backslash U)<m$ such that
$U$ is a Nash submanifold and $f$ is a Nash function on $U$. 
}\smallskip 

\begin{proof} Let 
 $$\Gamma(f)\subset S\times \RR$$
 be the graph of $f$, which is a closed semi-algebraic set;
 the projection $p_1: \Gamma(f)\to S$ is a semi-algebraic homeomorphism.
 
 By (\PropIsopartdense),  there is an open semi-algebraic set $U$ of $S$ 
 with $\dim (S\backslash U)<m$, which is a Nash submanifold, such that 
$p_1^{-1}(U) $ is also a Nash submanifold and 
 $p_1: p_1^{-1}(U) \to U$ is a Nash diffeomorphism.
 If $f_0=f|_U : U\to \RR$ is the restriction of $f$ to $U$, then
 $p_1^{-1}(U)$ is identified with the graph of $f_0$, so it follows that
 $f_0$ is a Nash map.
\end{proof}

\noindent (\DefIntegralSalgset)\quad
Let $S$ be a semi-algebraic set of dimension $m$.
We shall define the space of {\it generically defined} smooth differential forms  on $S$.
Denote by $Op(S)$ the collection of non-empty semi-algebraic subsets $U$ of $S_{reg}$ with $\dim (S\minus U)<m$, which is a Nash submanifold of dimension
$m$.
Note that if $U, U'\in Op(S)$, then 
$U\cap U'\in Op(S)$. 
Thus with respect to the ordering $U\le U'\Leftrightarrow U\supset U'$, 
$Op(S)$ is a directed set.

For $U$ in $Op(S)$, $\cA^p(U)$ denotes the vector space of smooth $p$-forms on  $U$. 
If $U'\subset U$ are subsets in $Op(S)$, there is the
restriction map $\cA^p(U) \to \cA^p(U')$.
Set  
$$A_S^p=\varinjlim \cA^p(U)\,,$$
the inductive limit over $U\in Op(S)$.
There is the differential map $d: A_S^p\to A_S^{p+1}$ induced from 
the differential $d: \cA^p(U)\to \cA^{p+1}(U)$.
There is also the wedge product map $A_S^p\otimes A_S^q\to A_S^{p+q}$.

Suppose that an orientation of $S_{reg}$ is given
(we will say for simplicity that $S$ is oriented).
If we are given an element $\vphi$ of $A_S^m$
represented by $\vphi\in \cA^m(U)$ for $U\in Op(S)$, 
and a Borel measurable function $f$ on a Borel measurable set $E$
of $S$,  then the restriction $f|_{U\cap E}$ is a Borel measurable
function on $E\cap U$.  We can thus apply the definition of (\Notationintegral)
to the Nash submanifold $U$, the form $\vphi$ and 
the restriction $f|_{U\cap E}$, so we have the integral
$\int_{E\cap U} (f|_{E\cap U})d\mu_\vphi$. By (\Propmeasurezero), 
this is well-defined, independent of the choice of a representative for $\vphi\in A^m_S$.
Thus we write
$$\int_E f d\mu_\vphi \,\quad \text{for}\quad \int_{E\cap U} (f|_{E\cap U}) d\mu_\vphi\,.$$
If $f=1$, it will be abbreviated to $\int_E \vphi$. 

Similarly, for $S_{reg}$ oriented or not, we have the integral $\int_{E\cap U} (f|_{E\cap U})
 d\mu_{|\vphi|}$
 and for which we will write
$\int_E f d\mu_{|\vphi|} $. 

In the above discussion we could have replaced ``smooth" forms by ``measurable" forms on $U$; but we will be mainly interested in the former.
\bigskip 

\noindent (\Integralinvariance)\quad
The space $A_S^m$ and the integral is a topological invariant in the following sense.
Let $f: S\to S'$ be a semi-algebraic homeomorphism of semi-algebraic sets of dimension
$m$.   There exists an isomorphism $f^*: A_{S'}^m\to A_S^m$ given as follows.
For $\vphi\in \cA^m(U')$ with $U'\in Op(S')$, we may assume by (\PropIsopartdense)
that $U=f^{-1}(U')$ is in $Op(S)$ and $U\to U'$ is a Nash diffeomorphism; then one has
$(f|_U)^*\vphi\in \cA^m(U)$, and represents an element $f^*\vphi$ of $A^m_S$.

Assume further that $S$ and $S'$ are oriented, and that $f$ preserves the orientation
(on suitable $U$ and $U'$). Then from the definitions we have
$$\int_S f^*\vphi= \int_{S'} \vphi\,.$$
Similarly in the unoriented case, $\int_S |f^*\vphi|= \int_{S'} |\vphi|$. 
\bigskip

\newpage

{\bf \S \sectResult. Basic results on integration of semi-algebraic sets}\bigskip

\noindent (\Pullback)\quad
 Let $S\subset \RR^n$ be a semi-algebraic set of dimension $\le m$.
Let $h: S\to \RR^\ell$ be a continuous semi-algebraic map, 
$\Omega\subset \RR^\ell$ be an open set containing $h(S)$,
 and $\psi$ be a smooth
$m$-form on $\Omega$. 

By (\CorPropIsopartdense) there exists a set $U\in Op(S)$ such that 
the restriction $h|_U :U\to \RR^\ell$ is a Nash map. 
The map $h|_U: U \to \Omega$ is a smooth map, so we can define the pull-back $(h|_U)^*\psi$ as a smooth form on $U$, namely as an element of $\cA^m(U)$.   We denote its image in $A^m_S$ by $h^*\psi$. 
(If $\dim S<m$, we have of course $h^*\psi=0$.)
 When $\dim S=m$,  by the previous section we have the integral $\int_S |h^*\psi|$, and if $S$ is oriented,  $\int_S h^*\psi$.  When $\dim S< m$, we set $\int_S |h^*\psi|=0$
 for convenience.

If $T=h(S)\subset \RR^l$ is the image of $h$, the map $h$ factors as 
$S\mapr{h} T\mapr{i} \RR^\ell$, where $i$ is inclusion.
Since $T$ has dimension $\le m$, $i$ is continuous semi-algebraic and $\Omega$ contains $T$, 
the pull-back $i^*\psi$ 
(also denoted $\psi|_T$) in  $A_T^m$ can be 
defined.  
\bigskip 

\noindent (\Defdelta) {\bf Definition.}\quad
Let $g: S\to T$ be a continuous semi-algebraic map between semi-algebraic sets.
Let 
$T_0=\{u\in T\, |\, g^{-1}(u) \mbox{ is a finite set}\}$ and 
$$\delta(g):= \operatorname{Max}\limits_{u\in T_0}\, \sharp(g^{-1}(u)) $$
which is a non-negative integer.
If  $T_0$ is empty (namely if all the fibers have positive dimension),  we set the number $\delta(g)$ to be zero.
\bigskip 

\noindent (\LemmaPartition) {\bf Lemma.}\quad{\it Let $S$ be a semi-algebraic set of dimension $m$, 
$\{S_i\}$ be a finite partition of $S$ into semi-algebraic subsets, 
and assume given, for each $S_i$ with dimension $m$, a set $U_i\in Op(S_i)$.
Then  
there exists a set $U\in Op(S)$ such that 
$U\cap S_i\subset U_i$ for each $S_i$ with dimension $m$, 
and such that 
$U$ is the disjoint union of $U\cap S_i$.
}

\begin{proof}
Let $V=S_{Zar}$ and $V_i=(S_i)_{Zar}$ be the Zariski closures of $S$ and $S_i$, respectively;
let $I$ be the subset of those indices $i$ with $ \dim S_i=m$, and 
$$V':=\bigcup_{i, j\in I} (V_i\cap V_j) \cup \bigcup_{i\not\in I} V_i\,.$$
For distinct $i, j\in I$, $\dim (V_i\cap V_j)<m$, hence $\dim V'<m$. 
Note that the sets $V_i-V'$ for $i\in I$ are disjoint from each other.
Thus the set
$$U:= \bigcup_{i\in I} (U_i-V')$$
satisfies the required property.
\end{proof}

\noindent (\PropPartitionforms) {\bf Proposition.}\quad{\it
Let $S\subset \RR^n$ be a semi-algebraic set of dimension $\le m$, 
and $h: S\to \RR^\ell$, $\Omega\subset \RR^\ell$, $\psi\in \cA^m(\Omega)$ be as above.
For $\{S_i\}$ a finite partition of $S$ into semi-algebraic subsets,
let $h_i: S_i\to \RR^l$ be the restriction of $h$ to $S_i$. 
Then we have
$$\int_S |h^*\psi| =\sum \int_{S_i} |h_i^*\psi|\,.$$
}

\begin{proof}
Let $U_i\in Op(S_i)$ be such that $h_i$ is a Nash map from $U_i$ to $\RR^\ell$. 
One takes $U\in Op(S)$ satisfying the conditions of the previous lemma.
By definition we have
$$\int_S |h^*\psi| =\int_U |(h|_U)^*\psi|\,\quad\text{and}\quad
\int_{S_i} |h_i^*\psi| =\int_{U_i} |(h_i|_{U_i})^*\psi|\,.$$
Since $U$ is the disjoint union of $U\cap S_i$, we have
$$\int_U |(h|_U)^*\psi|= \sum_i \int_{U\cap S_i} |(h|_U)^*\psi|\,.$$
For each $i$, 
$$\int_{U\cap S_i} |(h|_U)^*\psi| =\int_{U_i} |(h|_U)^*\psi|$$
by (\Propmeasurezero), 
proving the proposition.
\end{proof}

\noindent (\Propboundintegral) 
{\bf Proposition.}\quad{\it Under the assumption (\Pullback), we have
$$\int_S |h^*\psi| \le \delta (h) \int_T |(\psi|_T) |\,.$$
}\smallskip 

\begin{proof} 
(i) Let $\{T_i\}$ be a finite partition of $T$ into semi-algebraic subsets and 
$$h_i: S_i=p^{-1}(T_i) \to T_i\injto \RR^l$$
be the restriction of $h$ to $S_i$.  Then we claim that
the inequalities as stated for $h_i$'s imply the inequality for $h$.

Indeed, by (\PropPartitionforms)  we have
$$\int_T |\psi|_T|=\sum_i \int_{T_i} |\psi|_{T_i}|
\quad\text{and}\quad \int_S |h^*\psi|=\int_{S_i} |h_i^*\psi|\,,
$$ and $\delta(h)=\operatorname{Max}_i \delta(h_i)$.

(ii) Let $f: S\to S'$ be a semi-algebraic homeomorphism and $h': S'\to \RR^\ell$
be the semi-algebraic map such that $h'\scirc f=h$. 
Then $h^*\psi= f^* {h'}^*\psi $, where $f^*$ was defined in (\Integralinvariance), 
hence 
$$\int_S |h^*\psi|=\int_{S'} |{h'}^*\psi|\,.$$

For the proof of the proposition, take a finite partition $\{T_i\}$ of $T$ into semi-algebraic 
subsets
and semi-algebraic homeomorphisms $h^{-1}(T_i)\isoto T_i \times F_i$ over $T_i$ for some algebraic sets $F_i$ (see [BCR, p. 221, Theorem 9.3.2]).
By (i) we have only to prove the inequality for each $h_i: S_i\to T_i\subset \RR^\ell$. 

Thus assume that there exists a semi-algebraic homeomorphism
$S\to T\times F$ over $T$.
By (ii) we may assume that
$S=T\times F$ and the map $S\to T$ is projection $p: T\times F\to T$.

If $\dim T<m$, then $h^* \psi = p^*(i^*\psi)$ is zero since $i^*\psi=0$ for 
dimension reasons. Thus the integral of $|h^* \psi|$ on $S$ is zero, 
and the inequality holds.

If $\dim T=m$,  then $F$ consists of a finite number of points $\{P_j\}$, and
the restriction of $p$ to $U_j=T\times \{P_j\}$ is the identity map on $T$; thus
we have 
$$\int_S |h^* \psi |=\sum \int_{U_j} |h^* \psi |
=\sharp(F)\cdot \int_T |(\psi|_T)|\,,$$
so the stated inequality (indeed equality) holds.
\end{proof}

When $S$ is compact, the integral $\int_S h^*\psi$ absolutely converges:
\bigskip 

\noindent (\ThmNashfinite) 
{\bf Theorem.}\quad{\it Under the assumption  (\Pullback), assume that 
$S$ is compact. 
 Then the integral
$$\int_S |h^*\psi|$$
is absolutely convergent. 
(Thus if $S$ is oriented, $\int_S h^*\psi\in \RR$ is defined.)
In particular, if $\psi$ is a smooth $m$-form on an open neighborhood of 
$S$ in its ambient space $\RR^n$, then $\int_S |(\psi|_S)|$
is absolutely convergent. 
}\smallskip 

\begin{proof} Let $1\times h: S\injto \RR^n\times \RR^\ell$ be the semi-algebraic 
inclusion obtained as the product of the inclusion $S\to \RR^n$ and the map $h$.
Then $h^*\psi =(1\times h)^* p_2^* \psi$, where $p_2: \RR^n\times \RR^\ell\to R^\ell$
is the second projection. 
Replacing $S\subset \RR^n$ with $S\subset \RR^n\times \RR^\ell$,  we will 
assume that $h$ is the inclusion to its ambient space, and $\psi$ is defined on an
neighborhood of $S$.

The form $\psi$ is a finite sum of the forms
$$ a \,dx_{i_1}\wedge\cdots\wedge dx_{i_m}$$
where $a$ is a smooth function defined on a neighborhood of $S$.
By (2) of (\Propsumphi), one may assume $\psi=a \,dx_{i_1}\wedge\cdots\wedge dx_{i_m}$. 
Since $a$ is bounded on $S$ by compactness,  it is enough to 
consider the case
$\psi=dx_1\wedge\cdots\wedge dx_m$.

Let 
$q: S\to \RR^m$
be composition of the inclusion with projection $\RR^n\to \RR^m$
to the first $m$ coordinates, and apply the proposition
to $q$ and $\psi=dx_1\wedge\cdots\wedge dx_m$. 
We obtain
$$\int_S |(\psi|_S)| \le \delta (q) \int_{q(S)} | dx_1\wedge\cdots\wedge dx_m |\,,$$
and the integral on the right hand side is finite since $q(S)$ is compact so
its Lebesgue measure in $\RR^m$ is finite.
\end{proof}

\noindent (\ThmNashfinite.1) {\bf Remark.} \quad 
An analogous statement for a smooth manifold $S$ with boundary is false.
More precisely, let $S$ be a compact $m$-dimensional manifold with boundary, and let $S\injto \RR^n$ be
a continuous embedding which is smooth on the interior $\ctop{S}$ and on the boundary $\delta S$
of $S$. For $\psi$ a smooth $m$-form on $\RR^n$, the 
integral $$\int_{\ctop{S}} |(\psi|_{\ctop{S}})\,|$$ can be infinite.
\bigskip

For the proof of the Cauchy formula, and for the proof of 
the next proposition, it 
 is convenient to state a consequence of (\Propboundintegral) as
 in  the following proposition.   
 Note that, if $f$ is a continuous semi-algebraic function, then $df\in A_S^1$.
\bigskip 

\noindent (\PropBoundingAdd)
{\bf Proposition.}\quad{\it
Let $S$  be a compact oriented semi-algebraic set of dimension $m$, and 
$a$, $f_1, \cdots, f_m$ be continuous semi-algebraic functions on $S$, 
so $a \, df_1\wedge\cdots\wedge df_m\in A_S^m$.
We take  $f=(f_1, \cdots, f_m): S\to \RR^m$, a continuous semi-algebraic
map.
Then we have
$$\left|\int_S a\, df_1\wedge\cdots\wedge df_m\right|\le 
\delta(f) \cdot (\Max_{x\in S}|a(x)|)\cdot \int_{f(S)} dx_1\cdots dx_m\,.$$
}\smallskip 

\begin{proof}
We obviously have
$$\left|\int_S a\, df_1\wedge\cdots\wedge df_m\right|\le
\Max_{x\in S}|a(x)|\cdot \int_S | df_1\wedge\cdots\wedge df_m |\,.$$
Applying (\Propboundintegral) to the map $f$ and the form $\psi=dx_1\wedge\cdots\wedge
dx_m$, we have
$$\int_S | df_1\wedge\cdots\wedge df_m |\le
 \delta(f)\cdot\int_{f(S)}| dx_1\wedge\cdots\wedge dx_m |\,,$$
verifying the assertion. 
\end{proof}

\noindent (\Propcontinuity)
{\bf Proposition.}\quad{\it
Let $S\subset \RR^m$ be a compact oriented semi-algebraic set of dimension $m$,  
$ H: S\times [0, 1) \to  \RR^\ell$
be a continuous semi-algebraic map, $\Omega$ an open neighborhood of 
the image of $H$, and $\psi$ a smooth $m$-form on $\Omega$. 
 For $t\ge 0$, denote by 
$h_t: S\to  \RR^\ell$ the restriction of $h$ to $S\times \{t\}$.
Then we have
$$\int_S h_t^*\psi \to \int_S h_0^*\psi \quad\text{as $t\to 0$}\,.$$
}\smallskip 

\begin{proof}
By Theorem (\ThmNashfinite), we know that for each $t$,  $\int_S h_t^*\psi$ is absolutely convergent.

One may assume that the image of $H$ is contained in $\Omega$, and 
$\psi=c\, dx_{i_1}\wedge\cdots\wedge dx_{i_m}$ with $c\in \cA^0(\Omega)$;
after renumbering, let $\psi=c\, dx_1\wedge\cdots\wedge dx_m$. 
The functions $f_i(x, t)=H^*(x_i)$ for $i=1, \cdots, m$ and $a(x, t)=H^*c$ are 
continuous semi-algebraic on $S\times [0, 1)$.
Our assertion is now
$$\int_S a(x, t)\, df_1(x, t)\wedge\cdots\wedge df_1(x, t) 
\to \int_S a(x, 0)\, df_1(x, 0)\wedge\cdots\wedge df_1(x, 0) 
$$
as $t\to 0$.  
Set 
$I(a; f_1, \cdots, f_m)=\int_S a\, df_1\wedge\cdots \wedge df_m$.
We have 
$$\begin{array}{cl}
&I(a(x, t); f_1(x, t), \cdots, f_m(x, t))-I(a(x, 0); f_1(x, 0), \cdots, f_m(x, 0)) \\
=&I(a(x, t)-a(x, 0); f_1(x, t), \cdots, f_m(x, t))\\
&+I(a(x, 0); f_1(x, t)-f_1(x, 0), \cdots, f_m(x, t)) \\
&+I(a(x, 0);f_1(x, 0),  f_2(x, t)-f_2(x, 0), f_3(x, t)\cdots, f_m(x, t)) \\
&\quad\vdots  \\
&+I(a(x, 0);f_1(x, 0),  \cdots, f_{m-1}(x, 0), f_m(x, t)-f_m(x, 0))\,.
\end{array}
$$
We will examine each term in the sum.

Let 
$$H_0=(f_1(x, t), \cdots, f_m(x, t)): S\times [0,1) \to \RR^m\times [0, 1)$$
and $h_{0, t}=H_0|_{S\times \{t\}}: S\to \RR^m$. 
Note that $\delta(h_{0, t})\le \delta(H_0)$. 
Also we have a map, for $i=1, \cdots, m$, 
$$H_i=(f_1(x, 0), \cdots, f_{i-1}(x, 0), f_i(x, t)- f_i(x, 0), f_{i+1}(x, t), \cdots, 
f_m(x, t)): S\times [0,1) \to \RR^m\times [0, 1)$$
and its restriction $h_{i, t}=H_i|_{S\times \{t\}}: S\to \RR^m$;
we have $\delta(h_{i, t})\le \delta(H_i)$.

For the first term in the sum,  we have by the previous proposition
$$\begin{array}{cl}
&|I(a(x, t)-a(x, 0); f_1(x, t), \cdots, f_m(x, t))| \\
&\le \delta(h_{0, t})\cdot\Vert a(x, t)-a(x, 0)\Vert_S \cdot vol (h_{0, t}(S))\,,
\end{array}$$
where $\Vert a\Vert$ denotes the sup norm of a function on $S$,
and $vol$ is the Lebesgue measure in $\RR^m$.
One has $\delta(h_{0, t})\le \delta(H_0)$, $\Vert a(x, t)\Vert_S\to 0$ 
as $S$ is compact, and $vol (h_{0, t}(S))$ is bounded by $vol (H_0(S\times [0, 1/2])\,)$
for $t\le 1/2$; thus the first term tends to zero.

For the other terms, we have 
$$\begin{array}{cl}
&|I(a(x, 0);f_1(x, 0), \cdots, f_{i-1}(x, 0), f_i(x, t)- f_i(x, 0), f_{i+1}(x, t), \cdots, 
f_m(x, t) \,)| \\
&\le \delta(h_{i, t})\cdot\Vert a(x, 0)\Vert_S \cdot vol (h_{i, t}(S))\,.
\end{array}$$
We have $\delta(h_{i, t})\le \delta(H_i)$.
Since $f_i(x, t)\to f_i(x, 0)$ as $t\to 0$  uniformly on $S$ and $f_j(x, t)$ is uniformly bounded for each $j$, $vol (h_{i, t}(S))\to 0$. 
Hence the terms also tend to zero.
\end{proof}

Recall that $\Delta^m$ is a compact semi-algebraic set of $\RR^m$ 
defined by $x_i\ge 0$ for $i=1, \cdots, m$ and $\sum x_i\le 1$.
We equip $\Delta^m$ with the same orientation as $\RR^m$. 
We say that a continuous semi-algebraic map $h: \Delta^m\to \RR^n$
is {\it facewise smooth} if the restriction of $h$ to the interior of any
face is a smooth map.

It can be proven that any compact semi-algebraic 
set $S\subset \RR^n$ allows a semi-algebraic triangulation such that the inclusion of each simplex
 in $\RR^n$ is facewise smooth.
\bigskip 

\noindent (\Thmstokes)
{\bf Theorem.}\quad 
Let $S=\Delta^m$ and $h: S \to \RR^\ell$ be a continuous semi-algebraic, facewise
smooth map; let $\psi$ be a smooth $(m-1)$-form defined on an open 
neighborhood of $h(S)$ in $\RR^\ell$. 
Denote by $i: \delta S \injto S$ the embedding of the boundary of $S$.
Then we have the identity (the Stokes formula):
$$\int_{\delta S}(i\scirc h)^* \psi =\int_{ S}  h^*(d\psi)\,.$$
\smallskip 

\begin{proof} 
For $0\le t <\rho=1/(m+1)$, let 
$$S_t=\{(x_1, \cdots, x_m) \in S|\,  x_i\ge t\, \text{ and}\,  \sum x_i \le 1-t\}\,$$
be the ``$t$-excision" of $S$. 
Let 
$$r=(r_1(x, t), \cdots, r_m(x, t)): S\times [0, \rho) \to S$$
be the map given by $r_i(x, t)= (1-(m+1) t)x_i +t$. 
Then $r$ is a continuous semi-algebraic map,   
the restriction $r_t=r|S\times \{t\}$ is a semi-algebraic 
homeomorphism from $S$ to $S_t$, and $r_0$ is the identity map.

For $t>0$, the embedding $S_t\injto \RR^n$ is a smooth map,
so $h^*\psi$ is a smooth form on $S_t$, a manifold with corners.
(This is where we use the assumption of facewise smoothness.
So we only need that the map $i$ be smooth on the interior of 
$S$.) 
By the Stokes formula for a smooth form on a manifold with corners, we
have
$$\int_{\delta (S_t)} (hi)^*\psi
=\int_{ S_t} h^*(d\psi)\,.\eqno{(*)}$$

Since $r_t: \delta S\to \delta (S_t)$ is a diffeomorphism of manifolds, 
the left hand side of (*) equals
$$\int_{\delta S} r_t^* (hi)^*\psi\,.$$
Consider the composition of the maps
$H: \delta S\times [0, \rho)\mapr{r}\delta S \mapr{hi} \RR^\ell$;
applying (\Propcontinuity) we obtain
$$\int_{\delta S} r_t^* (hi)^*\psi \to \int_{\delta S}  (hi)^*\psi\,.$$

As for the right hand side of (*), since $S_t$ is increasing and 
$\cup S_t =\ctop{S}$, by Lebesgue's monotone convergence theorem we have
$$ \int_{ S_t} h^*(d\psi) \to \int_{\ctop{S}} h^*(d\psi)$$
as $t\to 0$, and $\int_{\ctop{S}} h^*(d\psi)=\int_{S}h^*(d\psi)$ since 
the complement of $\ctop{S}$ has dimension $<m$.
\end{proof}
\bigskip

\newpage

{\bf \S \sectAllow. Integrals on allowable semi-algebraic sets}
\bigskip 

We begin by introducing some notation.
For integers $0\le p\le n$,
let $\RR^n=\RR^p\times \RR^{n-p}$ be Euclidean space with coordinates $(r_1, \cdots, r_p,  x_{p+1}, \cdots, x_n)$. 
The first $p$ coordinates and the last $n-p$ coordinates will play 
 different roles. 
Let $H_i=\{r_i=0\}$ for $1\le i\le p$.  A {\it face} is the intersection of 
some of the $H$'s; thus a face is of the form 
$H_I=\cap_{i\in I} H_i $ for a subset $I$ of $\{1, \cdots, p\}$.  We include $\RR^n$
as a face. 
A face $H_I$ has induced coordinates, given by $r_i$, $i\not\in I$,
 and all of $x_{p+1}, \cdots, x_n$;
 thus $H_I=\RR^k\times \RR^{n-p}$ if $k=p-|I|$. 
Note that $H_1\cap \cdots\cap H_p$ is the smallest face, and its dimension is $n-p$.  

For $m\le n$, by an $m$-form with logarithmic singularities 
(or just a logarithmic $m$-form) on $\RR^n$,
we mean a smooth $m$-form defined on $\RR^n-(\cup H_i)$ 
which  can be written 
$$\omega=
\sum_{I=(i_1, \cdots, i_m)} a_I(r, x) 
\frac{dr_{i_1}}{r_{i_1}}\wedge\cdots\wedge \frac{dr_{i_k}}{r_{i_k}}\wedge dx_{i_{k+1}}\wedge\cdots\wedge dx_{i_m}\,$$
where the sum is over the sequences 
$1\le i_1<\cdots <i_k\le p< i_{k+1}<\cdots i_m\le n$ and 
$a_I(r, x)$ are smooth functions  
in $(r, x)\in \RR^n$.

Let  $A$ be a closed semi-algebraic set of $\RR^n$ 
such that  $\dim A=n$ and $\dim (A\cap(\cup H_i)\,)<n$. 
For a logarithmic $n$-form $\omega$ on $\RR^n$, we have defined in \S 1 
the integral
$$\int_A|\omega|\,,$$
also written $vol(A; \omega)$, 
which  is a non-negative real number (possibly $+\infty$).
We will study when it is finite (then we say that $\int_A \omega$ is 
absolutely convergent).
\bigskip

The line of argument leading to the main theorem of this section, (\Thmconvergence), is as follows.
We shall first the notion of allowability and (almost) strict allowability.
\smallskip 

(1)  We show Theorem (\thmreduction), 
which states that an allowable semi-algebraic set can be made, by a succession of ``permissible" blow-ups, almost strictly allowable. 

(2) Given an almost strictly allowable compact semi-algebraic set $A$ in $\RR^n=\RR^p\times \RR^{n-p}$, one shows:
after an appropriate linear change of variables in $(x_{p+1}, \cdots, x_n)$, 
the projection map $pr_1:\RR^n\to \RR^{n-1}$,
$(r_1, \cdots, r_p, x_{p+1}, \cdots, x_n)\mapsto (r_1, \cdots, r_p, x_{p+1}, \cdots, x_{n-1})$,
when restricted to $A$, 
has finite fibers.
See (\propfiniteness.1).

(3) When this finiteness property for projection $pr_1$ holds for $A$, one can prove the absolute 
convergence of $\int_A \omega$. See (\PropLtypeineq) and the proof of (\Thmconvergence.1).

(4) Let $A$ be an allowable compact semi-algebraic set in $\RR^p\times \RR^{n-p}$. Combining (1)-(3), it 
follows that the integral  $\int_A \omega$ is absolutely convergent.  This is Theorem (\Thmconvergence).
\bigskip  

Let $C\subset \RR^n$ be a  closed semi-algebraic subset. 
We denote by $C_{Zar}$ the Zariski closure 
of $C$.
The dimension of $C$ is, by definition, the dimension of its Zariski closure: 
$\dim C= \dim C_{Zar}$. 
\bigskip 

\noindent (\Defallowable) {\bf Definition.}\quad 
Let $C$ be a  closed semi-algebraic subset of $\RR^n$. 
We say $C$ is {\it allowable} with respect to a face $F$  
of $\RR^n$ if $(C\cap F)_{Zar}\subsetneq F$
(equivalently if  $\dim (C\cap F)< \dim F$). 

We define $C$ to be allowable in $\RR^n$  if it is allowable with respect to any proper face of $\RR^n$. 
\bigskip 

\noindent (\DefAllowable) {\bf Definition.}\quad
Let $C$ be a  closed semi-algebraic subset of $\RR^n$.
We say $C$ is {\it strictly allowable}  with respect to a face $F$ 
if $(C\cap F)_{Zar}$ does not contain any face.

The set $C$ is defined to be {\it almost strictly allowable}
(resp. {\it strictly allowable}) in $\RR^n$ if $C$ is strictly allowable with respect to any proper face (resp. any face including $F=\RR^n$).

Note that if $C$ is strictly allowable in $\RR^n$, then 
$\dim C<n$. 
But a set $C$ with $\dim C=n$ may be allowable or almost 
strictly allowable. 
\bigskip

\noindent (\basicproperties) {\bf Basic properties.}\quad 
(1) 
It is clear that if $C$ is strictly allowable with respect to a face $F$,  
then it is allowable with respect to $F$. 

Also, if $C$ is 
strictly allowable in $\RR^n$, then it is   almost 
strictly allowable in $\RR^n$; if $C$ is almost strictly allowable in $\RR^n$, then it is allowable  in $\RR^n$.

(2) If $C$ is strictly allowable with respect to $F$, and 
$F'$ is a face contained in $F$, then $C$ is also 
strictly allowable with respect to $F'$.
(This follows from the inclusion $(C\cap F)_{Zar}\cap F'\supset (C\cap F')_{Zar}$. )
Thus $C$ is almost strictly allowable in $\RR^n$ if it is strictly allowable with respect to any codimension one face of $\RR^n$, 
and it is strictly allowable in $\RR^n$ if it is strictly allowable with respect to $F=\RR^n$. 

An analogous statement for allowability is false: if $C$ is allowable with respect to a face, it may not be allowable 
with respect to a smaller face. 

(3) From the definition,  $C$ is allowable (resp. strictly allowable) with respect to $F$ if and only if $C\cap F$ is allowable (resp. strictly allowable) with respect to $F$. 

Thus, if $C$ is contained in a proper face $F$, $C$ is strictly allowable
 in $\RR^n$ if and only if $C$ is strictly allowable in $F$. 
Because of this we often do not mention the ambient space of $C$
 for the condition 
of strict allowability. 
(The corresponding statement does not hold if we replace strict allowability
with allowability or almost strict allowability.)

(4) If $C'$ is a  closed semi-algebraic subset of $C\subset \RR^n$, 
allowability or (almost) strict allowability for $C$ obviously implies 
the same for $C'$. 

If $C$ is the union of  closed semi-algebraic subsets $C_1$ and $C_2$, then 
allowability or (almost) strict allowability for $C_1$ and $C_2$ implies 
the same for $C$.

(5) Assume $C$ is an algebraic subset of $\RR^n$. 
We say that $C$ {\it meets the faces properly} if for  each face $F$, $\dim (C\cap F) \le \dim C+ \dim F-n$.  
If $C$ meets the faces properly and if $\dim C<n$, 
then $C$ is strictly allowable.  
\bigskip 

We will also need to consider semi-algebraic sets of blow-ups of $\RR^n$, 
as well as the extended notion of allowability.
Denote by $S$ the Euclidean space  $\RR^n$ as above. 
By a {\it permissible} blow-up of $S$ we mean the blow-up of 
a face.  Thus it is of the form 
$\mu=Bl_\Sigma: \tilde{S}=\tilde{\RR}^n\to \RR^n$, 
where $\Sigma$ is a  face of codimension $\ge 2$. 
On $\tilde{S}$ we have the strict transforms $H'_i$ of $H_i$ for $0\le i\le p$, 
and the exceptional divisor $E$, which form a normal crossing divisor. 
Thus on $\tilde{S}$ again one has the notion of faces (note the minimal faces are 
of dimension $n-p$)
and permissible blow-ups.  
Iterating this we have a succession of permissible blow-ups
$$\mu=\mu_1\mu_2\cdots\mu_N:
 \ti{S}=S_N\mapr{\mu_N}S_{N-1}\to\cdots
\to S_1 \mapr{\mu_1}S_0=S$$
where $\mu_{k+1}$ is the blow-up of a face $\Sigma_k$ of $S_{k-1}$. 
The $\ti S$ is a real algebraic variety (in the sense of 
[BCR, Definition 3.2.11]).
On $\tilde{S}$ one has the notion of faces; the minimal faces are 
of dimension $n-p$.  
There is a coordinate 
chart containing each minimal face. 

Since $\ti S$ is a real algebraic variety, there is the notion
of semi-algebraic sets of $\ti S$ (see [BCR, p. 64]);
note that a set is semi-algebraic if and only if its intersection with 
each chart is semi-algebraic.
For a closed semi-algebraic set $C\subset \ti{S}$, one can define the notion of 
(strict) allowability with respect to a face and (strict or almost strict) allowability in  $\ti{S}$
in the same manner as for $\RR^n$.
The basic properties above hold true without any change. 
We add another property to the list:
\bigskip

(6) A set $C$ of $\ti S$ is  allowable (resp. strictly 
allowable, resp. almost strictly allowable) 
if and only if for each chart $U$, $C\cap U$ is allowable (resp. strictly allowable, resp. almost strictly allowable) in $U$. 
\smallskip

Indeed, for a face $F$, 
$C\cap F$ is the union $\cup(C\cap U\cap F)$ for 
the charts $U$ of $\ti{S}$, thus
$$\dim (C\cap F)=\sup \dim (C\cap U\cap F)\,.$$
Hence follows the assertion for allowability.

To show the assertion for strict allowability, 
note that
$$(C\cap U\cap F)_{Zar}=(C\cap F)_{Zar} \cap U\,,$$
where $Zar$ on the left means the Zariski closure in
$U$.  It follows that if $(C\cap F)_{Zar}$ contains 
a face $\Sigma$, then for a chart $U$ such that $U\cap \Sigma$ is non-empty, $(C\cap U\cap F)_{Zar}$ contains 
the face $U\cap \Sigma$. Conversely, if $(C\cap U\cap F)_{Zar}$ contains a face of the form $U\cap \Sigma$, 
then $(C\cap F)_{Zar} \supset U\cap \Sigma$, and taking 
the Zariski closure in $F$ gives
$(C\cap F)_{Zar} \supset \Sigma$.
\bigskip

\noindent (\lemmabasic) {\bf Lemma.}\quad{\it 
Let $\mu: \tilde{\RR}^n\to \RR^n$ be a succession of permissible
blow-ups.  Let $G$ be a face of $\tilde{\RR}^n$ and let 
$F=\mu (G)$. 

(1) For a  closed semi-algebraic set $C\subset \RR^n$, if 
$(\mu^{-1}(C)\cap G)_{Zar}$ contains a face $\Sigma$, 
then  $(C\cap F)_{Zar}\supset \mu(\Sigma)$. 

(2) If $C$ is strictly allowable with respect to $F$, 
then $\mu^{-1}(C)$ is strictly allowable with respect to $G$.

(3) If $C$ is strictly allowable (resp. almost strictly allowable) in ${\RR}^n$, 
then $\mu^{-1}(C)$ is  strictly allowable (resp. almost strictly allowable)
in $\tilde{\RR}^n$.
}\smallskip 

\begin{proof}
One has an obvious inclusion 
$$(\mu^{-1}(C)\cap G)_{Zar}\subset \mu^{-1}((C\cap F)_{Zar})\cap G\,.$$
If  $(\mu^{-1}(C)\cap G)_{Zar}$ contains a face $\Sigma$, 
it follows that $\mu^{-1}((C\cap F)_{Zar})\supset \Sigma$, hence 
$(C\cap F)_{Zar}$ contains the face $\mu(\Sigma)$. 

The assertions (2) and (3) follow from (1). 
\end{proof}

Let $k$ be an infinite field and $S=\Spec k[x_1, \cdots, x_n]$ be affine $n$-space
over $k$.  Let $0\le p\le n$ be a given integer.  
For $i=1, \cdots, p$ let $H_i=\{x_i=0\}$ be
a coordinate hyperplane.  As in the real case we have the notion of faces 
and permissible blow-ups.
The following result is due to [Bl] (in case $p=n$). 
\bigskip 

\noindent (\blochthm) {\bf Theorem}\,[Bloch].\quad{\it  Let $V\subset S$ be a closed subvariety, not contained in 
any $H_i$ for $i=1, \cdots, p$. 
Then there is a permissible blow-up $\mu: \ti S\to S$ such that the strict transform 
of $V$ in $\ti S$ meets the faces properly.
}\smallskip 

\begin{proof} We indicate how to modify the proof of Theorem (2.1.2) in [Bl]
(which is the case $p=n$). 

If $\mu: \ti S\to S$ is a succession of permissible blow-ups, one has distinguished 
prime divisors on $\ti S$ that form a normal crossing divisor, and minimal faces 
are of dimension $n-p$.  Let $\cB_S$ be the full subcategory of $S$-schemes 
whose objects are permissible blow-ups $\mu: T\to S$. 
For each $T$ in $\cB_S$ we have the set of minimal faces $\cV(T)$ , and define
$\cV$ to be the inverse limit of $\cV(T)$ as $T$ varies over $\cB_S$,
$$\cV=\varprojlim_{T\in \cB_S} \cV(T)\,.$$
One has projection $pr_T: \cV\to \cV(T)$. 
With this one has:
\bigskip 

\noindent (\blochthm.1) {\bf Lemma.}\quad{\it For any $v\in \cV$, there is an object $\mu: T\to S$  in 
$\cB_S$ such that $pr_T(v)\not\subset V'$, where $V'$ is the strict transform of 
$V$ in $T$. }
\bigskip 

The proof is parallel to that for Lemma (2.1.2.1) in [Bl] except we modify 
the argument as follows. 
 
Let $V$ be a prime divisor given by $f(x_1, \cdots, x_n)=0$. 
We can view $f$ as a polynomial in  $x_1, \cdots, x_p$ with 
 coefficients polynomials in $k[x_{p+1}, \cdots, k_n]$. 
Then consider the set 
$$M=\{(r_1, \cdots, r_p)\in \ZZ_{\ge 0}^p\,|\, \mbox{the coefficient of $x^r$
of $f$ is non-zero.}\,\}.$$
Let 
$$\Delta\subset\RR_{\ge 0}^p:= \mbox{the convex hull of $\bigcup_{r\in M} (r+\RR_{\ge 0}^p)$}\,.$$
By the same argument as in {\it loc. cit.}, one shows that there is an object 
$\mu: T\to S$ such that the strict transform of $V$ has equation $f'$ having 
a non-zero constant term $c(x_{p+1}, \cdots, x_n)$; in other words $V'$ meets the faces properly. 
\end{proof}

A closed semi-algebraic set of $\RR^n=\RR^p\times\RR^{n-p}$ 
(or of a  permissible blow-up of $\RR^n$),
which is of dimension
$<n$ and almost strictly allowable,  can be made
strictly allowable by a succession of permissible blow-ups:
\bigskip 

\noindent (\propreduction) {\bf Proposition.}\quad{\it 
Let $S'\to \RR^n=\RR^p\times\RR^{n-p}$ be a 
succession of permissible blow-ups, and 
$B$ be a closed semi-algebraic subset of $S'$
with $\dim B <n$, which is
 almost strictly allowable in $S'$. 
Then there is a succession of permissible blow-ups $\mu:
\tilde{S}'\to S'$
 such that $\mu^{-1}(B)$ is strictly allowable.}
\smallskip 

\begin{proof}
{\it Step 1.}\quad We first consider the case 
$S'=\RR^p\times\RR^{n-p}$. 
If $B$ is contained in an $H_i$, then $B$ is strictly allowable by 
assumption, so the assertion trivially holds (with $\mu=id$). 
Thus we assume $B\not\subset H_i$ for each $i$. 

Let  $V=B_{Zar}$, which is an algebraic subset $\not\subset H_i$, and apply 
(\blochthm)(with $k=\RR$). 
There is a succession of permissible blow-ups 
$$\mu=\mu_1\mu_2\cdots\mu_N:
 \ti{\RR}^n=\cR_N\mapr{\mu_N}\cR_{N-1}\to\cdots
\to \cR_1 \mapr{\mu_1}\cR_0=\RR^n$$
where $\mu_{k+1}$ is the blow-up of a face $\Sigma_{k}$, 
such that the strict 
transform $V'$ of $V$ meets the faces properly; thus $V'$ is
 strictly allowable in $\ti{\RR}^n$ 
by (\basicproperties), (5). 

Let 
$B'$ is the closure in the Euclidean topology 
of $\mu^{-1}(B\cap U)$, where $U\subset \RR^n$ is the largest open set 
over which $\mu$ is an isomorphism. 
Since $B'\subset V'$, $B'$ is also strictly allowable. 

One has clearly 
$$\mu^{-1}(B)=B'\cup \bigcup_E(\mu^{-1}(B)\cap E)\eqno{(\propreduction.a)}$$
where $E$ varies over the exceptional divisors of $\mu$. 
We claim that each $\mu^{-1}(B)\cap E$ is strictly allowable. 
Indeed if $G$ is a face of $E$, then $\mu(G)$ is a proper face of $\RR^n$, thus 
$B$ is strictly allowable with respect to $\mu(G)$;
by (2) of (\lemmabasic) it follows that 
$\mu^{-1}(B)\cap E$ is strictly allowable with respect to $G$. 

Since each subset in the union (\propreduction.a) is strictly allowable,
$\mu^{-1}(B)$ is also strictly allowable in $\ti{\RR}^n$.

{\it Step 2.}\quad We consider the general case.  To each 
minimal face $c$ of $S'$ there corresponds a chart $U_c=\RR^p\times\RR^{n-p}$.  By Step 1, there is a succession 
of permissible blow-ups 
$\mu_c: \ti{U}_c\to U_c$ such that $\mu_c^{-1}(B\cap U_c)$ is 
strictly allowable in $\ti{U}_c$. 
Let $\mu_c: \ti{S}'_c \to S'$ be the succession of permissible
blow-ups obtained by blowing up the closures of the centers of the blow-ups in
$\ti{U}_c\to U_c$; thus one has $\mu_c^{-1}(U_c) =\ti{U}_c$. 

We take a succession of blow-ups $\mu: \ti{S}' \to S'$ that 
dominates all $\mu_c$ (see [Bl, Corollary (1.2.2)]). Then $\mu$
factors as $\ti{S}'\mapr{\mu'}\ti{S}'_c\to S'$, so we have 
a commutative diagram
$$\begin{array}{ccccc}
  \ti{S}'      &\mapr{\mu'}&\ti{S}'_c  &\mapr{\mu_c} &S'   \\
     \cup     &                &\cup    &                  & \cup       \\
\mu^{-1}(U_c) &\mapr{\mu'}&\ti{U}_c &\mapr{\mu_c} &\phantom{\quad .}U_c\quad .
  \end{array}
  $$
Since  $\mu_c^{-1}(B\cap U_c)$
is strictly allowable in $\ti{U}_c$, by (3) of (\lemmabasic), 
$\mu^{-1}(B)\cap \mu^{-1}(U_c)$ is strictly allowable in $\mu^{-1}(U_c)$.
Thus by (6) of (\basicproperties), $\mu^{-1}(B)$  is strictly allowable in any 
chart of $\mu^{-1}(U_c)$. 
This being the case for all $c$, $\mu^{-1}(B)$  is strictly allowable in any 
chart of $\ti{S}'$; by the ``if" part of (6) of (\basicproperties), $\mu^{-1}(B)$ is 
strictly allowable in $\ti{S}'$.
\end{proof} 

\noindent (\thmreduction) {\bf Theorem.}\quad{\it 
Let $A$ be a  closed semi-algebraic subset of $\RR^p\times \RR^{n-p}$ which is allowable. 
Then there is a succession of permissible blow-ups $\mu:
\tilde{\RR}^n\to \RR^n$ such that $\mu^{-1}(A)$ is almost strictly allowable. 
}\smallskip 

\begin{proof} 
We give a slightly generalized statement. 
Let $S' \to S=\RR^n$ be an object of $\cB_S$ (the notation as introduced in the proof of (\blochthm)\,)
and $A\subset S'$ be a closed allowable
semi-algebraic subset.
For an integer $d$ with $0\le d\le n-1$, consider the following claim. 
\smallskip 

$(\text{\bf Claim})_d$:\quad{\it For any such $S'$ and 
$A$,  
there is a succession of permissible blow-ups $\mu:
\tilde{S'}\to S'$ such that the inverse image 
$\mu^{-1}(A)$ is strictly allowable with respect to 
each face $G$ of $\tilde{S'}$ with dimension $\le d$.
}\smallskip 

For $S'=S=\RR^n$ and $d=n-1$, this is the statement of the theorem. 
If $d\le n-p$, the dimension of the minimal faces, then 
$(\text{Claim})_d$ obviously holds with $\mu=id$, 
since allowability and strict allowability are equivalent with respect
to a minimal face. 
We will show $(\text{Claim})_d$ by ascending induction on $d$. 

Assume $(\text{Claim})_{d-1}$ holds, so there is a succession of permissible
blow-ups $\mu:\tilde{S'}\to S'$ satisfying $(\text{Claim})_{d-1}$. 
Rewriting $S'$ for $\tilde{S'}$, and $A$ for $\mu^{-1}(A)$, 
one may assume that \smallskip

(\thmreduction.a)\quad $A$ is strictly allowable with respect to each face of dimension $\le d-1$.
\smallskip 

Assume that $A$ is not strictly allowable with respect to a face $F$ of dimension $d$.
The set $A\cap F$ in $F$ satisfies the assumption for 
Proposition (\propreduction).  Indeed $\dim (A\cap F)<\dim F$ by allowability
and $A\cap F$ is strictly allowable with respect to each face $F'\subsetneq F$
by (\thmreduction.a).
It follows that 
 there is a succession of blow-ups 
$$\nu_F: \ti{F}\to F$$ with centers $\Sigma_k$ in $F$ (or in its
strict transform) such that $\nu_F^{-1}(A\cap F)$ is strictly allowable. 
Let $$\nu:  \tilde{S}'\to S'$$
 be the succession of blow-ups of 
$S'$ with the same centers $\Sigma_k$ as for $\nu_F$. The restriction of $\nu$ to 
$\ti{F}$ coincides with $\nu_F$, so  $\nu^{-1}(A)\cap \ti F
$ is strictly 
allowable; thus $\nu^{-1}(A)$ is strictly allowable with respect to $\ti{F}$. 

If $G$ is a face of $\ti{S}'$ of dimension $d$ such that $\dim \nu(G)=d$ and 
$A$ is strict allowable with respect to $\nu(G)$, then by (2) of (\lemmabasic),
$\nu^{-1}(A)$ is strictly allowable with respect to $G$. 
Also,  $\nu^{-1}(A)$ is strictly allowable with respect to each
face $G$ such that  $\dim \nu(G)<d$, 
by assumption (\thmreduction.a) and (\lemmabasic).

Let $\mathfrak{n}(S'; A)$ denote the number of faces $F$ of dimension $d$ such that $A$ is not 
strictly allowable with respect to $F$. Then the above argument shows 
$$\mathfrak{n}(\ti{S}'; \nu^{-1}(A))<\mathfrak{n}(S'; A)\,.$$
Iterating this process, we find a succession of blow-ups $\mu: \hat{S}'\to S'$ with 
$\mathfrak{n}(\hat{S}'; \mu^{-1}(A))=0$, thus $\mu^{-1}(A)$ is strictly allowable with respect to all faces 
of dimension $\le d$. 
\end{proof}

Let $A$ be a closed semi-algebraic subset of 
$\RR^{p}\times \RR^{n-p}$ with $p<n$. 
For a point $(c_{p+1}, \cdots, c_{n-1})$ of 
$\RR^{n-p-1}$, consider 
the linear change of variables 
$$x_i\mapsto x'_i+ c_i x'_n \quad (i=p+1, \cdots, n-1), \qquad
x_n\mapsto x'_n\,.\leqno{(LC)}$$ 
Then $A$ is viewed a subset of $\RR^{p}\times \RR^{n-p}$ with 
variables $(r_1, \cdots, r_p, x'_{p+1}, \cdots, x'_n)$. 

Rewriting $(x_{p+1}, \cdots, x_n)$
for $(x'_{p+1}, \cdots, x'_n)$, 
let $pr_1: \RR^{n}\to \RR^{n-1}$ be the projection to the first $n-1$ 
coordinates $(r_1, \cdots, r_p, x_{p+1}, \cdots, x_{n-1})$, let $pr_2: \RR^{n-1}\to \RR^p$ be the projection to 
the first $p$ coordinates $(r_1, \cdots, r_p)$, and 
let  $pr=pr_2pr_1: \RR^{n}\to \RR^{p}$ be the composition.

For $A$ strictly allowable, 
we have the following finiteness result for the map $pr_1$ restricted to $A$. 
(For $p=0$ and $A$  an algebraic set, this is a well-known proposition in algebraic geometry.)
\bigskip 

\noindent (\propfiniteness) {\bf Proposition.}\quad{\it  
Let $A$ be a strictly allowable 
compact semi-algebraic subset of 
$\RR^{p}\times \RR^{n-p}$, $p<n$. 
Then there exists an open neighborhood $V$ of the origin in $\RR^p$ and a non-empty Zariski open set $\cW$ of $\RR^{n-p-1}$ satisfying 
the following property: 
\smallskip 

For any point $(c_{p+1}, \cdots, c_{n-1})$ 
of $\cW$, perform the change of variables 
indicated above. 
Let $pr_2^{-1}(V)\subset \RR^{n-1}$ and 
$A_V=A\cap pr^{-1}(V)\subset \RR^{n}$, so that one has 
a commutative diagram 
$$\begin{array}{ccc}
 A_V   &\injto  &\RR^n \\
\mapd{pr_1} &     &\mapdr{pr_1} \\
pr_2^{-1}(V)   &\injto &\RR^{n-1}  \\
\mapd{pr_2} &     &\mapdr{pr_2} \\
 V&\injto  &\phantom{\,.}\RR^p\,.
 \end{array}
$$  
Then the fibers of the induced map $pr_1: A_V\to pr_1^{-1}(V)$
are finite sets. 
}\smallskip 

\begin{proof}
Note first that the assertion is obvious if $A\cap H_1\cap\cdots \cap H_p$ is empty.  For then the image of $A$ by the projection $pr$ to 
$\RR^p$ does not contain the origin; taking a neighborhood $V$ of the origin disjoint from the image of $A$, the assertion obviously holds.

Thus we will assume $A\cap H_1\cap\cdots \cap H_p$ is non-empty.
By the strict allowability, 
$A_{Zar}$ does not contain the face $H_1\cap\cdots\cap H_p$. 
Since $A_{Zar}$ is the intersection 
of the zero locus of the polynomials $f(r_1, \cdots, r_p, x_{p+1}, \cdots, x_n)$ that 
vanish on $A$, there exists a polynomial 
$f$ such that $A\subset Z(f)$ and 
$f(0, \cdots, 0,x_{p+1}, \cdots, x_n)$ is not 
zero as a polynomial in $(x_{p+1}, \cdots, x_n)$. 

Write 
$$f(r_1, \cdots, r_p, x_{p+1}, \cdots, x_n)
=\sum_m a_m(r_1, \cdots, r_p)x^m$$
where $m=(m_{p+1}, \cdots, m_n)$ varies over 
multi-indices, $a_m(r_1, \cdots, r_p)$
are polynomials, and $x^m=x_{p+1}^{m_{p+1}}\cdots x_n^{m_n}$. 
Let $f_j(r, x)=\sum_{|m|=j} a_m(r_1, \cdots, r_p)x^m$,  the
homogeneous of degree $j$ part  with respect to $x$. 
One has $f(r, x)=\sum_j f_j(r, x)$.
Note that a linear change of coordinates (LC) 
is compatible with the decomposition $f=\sum f_j$ and also with 
the substitution $r_1=\cdots= r_p=0$.

Since $f(0, \cdots, 0, x_{p+1}, \cdots, x_n)
$ is not zero as a polynomial and vanishes on a non-empty set, 
 its homogeneous degree $d$ with respect to $x$ is positive. 
 (Note $d$ may well be strictly smaller than the homogeneous degree of 
$f(r, x)$ with respect to  $(x_{p+1}, \cdots, x_n)$.)
There is thus a non-empty Zariski open set $\cW$
of $\RR^{n-p-1}$ such that for any 
$(c_{p+1}, \cdots, c_{n-1})\in \cW$, 
the corresponding change of variables (LC) 
renders the coefficient of 
${x'_n}^d$  in 
$$f_d(0, \cdots, 0, x'_{p+1}+c_{p+1} x'_n, \cdots,
x'_{p+1}+c_{p+1} x'_n, \cdots, x'_n)$$
non-zero.  
We will rewrite $x_i$ for $x'_i$, and 
$f(r_1, \cdots, r_p, x_{p+1}, \cdots, x_n)$ 
for 
$$f(r_1, \cdots, r_p,x'_{p+1}+c_{p+1} x'_n, \cdots,
x'_{p+1}+c_{p+1} x'_n, \cdots, x'_n)\,.$$

For a point $T_0=(r_1^0, \cdots, r_p^0,  x_{p+1}^0, \cdots, x_{n-1}^0)
$ of $\RR^{n-1}$, 
let 
$$f(T_0, x_n)=f(r_1^0, \cdots, r_p^0,  x_{p+1}^0, \cdots, x_{n-1}^0, x_n)$$
be the polynomial in $x_n$ obtained by substituting the values of the  coordinates of $T$ for the first $(n-1)$ variables. 
If $T_0=(0, \cdots, 0, x_{p+1}^0, \cdots, x_{n-1}^0)
$ is any point of $\RR^{n-1}$ with the first $r$ coordinates zero, 
$$f(T, x_n)=f(0, \cdots, 0, x_{p+1}^0, \cdots, x_{n-1}^0, x_n)
=f_d(0, \cdots, 0, x_{p+1}^0, \cdots, x_{n-1}^0, x_n)$$
is a non-zero polynomial of degree $d$. 
There is thus an open neighborhood $U(T_0)$ of $T_0$ in $\RR^{n-1}$ such that 
for each point $T\in U(T_0)$, $f(T, x_n)$ is a non-zero polynomial in $x_n$
of degree $\ge d$. 

Now for each point $T_0\in pr_1(A)\cap pr_2^{-1}(\{0\})$, take its neighborhood
$U(T_0)$ as above. 
Since the set $pr_1(A)\cap pr_2^{-1}(\{0\})$ is compact, one can cover it 
by a finite number of such neighborhoods; let $U$ be the union of them. 
Then $U$ is an open set of $\RR^{n-1}$ containing $pr_1(A)\cap pr_2^{-1}(\{0\})$, such that for each $T\in U$, $f(T, x_n)$ is a non-zero polynomial in $x_n$. 

The set $pr_1(A)-U$ is compact and its image
$pr_2(pr_1(A)-U)$ is a closed set not containing the point $0$. 
So 
$$V:=\RR^p-pr_2(pr_1(A)-U)$$
is an open neighborhood of $0$, and 
$$pr_2^{-1}(V)\cap pr_1(A) \subset U\cap pr_1(A)\,.$$
Thus for $T\in pr_2^{-1}(V)\cap pr_1(A)$, $f(T, x_n)$ is a non-zero 
polynomial in $x_n$, so the intersection 
$Z(f)\cap pr_1^{-1}(T)$ is finite. 
\end{proof}

\noindent (\propfiniteness.1) {\bf Corollary.}\quad{\it 
Let $A$ be an {\em almost} strictly allowable 
compact semi-algebraic subset of 
$\RR^{p}\times \RR^{n-p}$. 
Then there exist an open neighborhood $V$ of the origin in $\RR^p$ and a non-empty Zariski open set $\cW$ of $\RR^{n-p-1}$ satisfying 
the following property: 
\smallskip 

 \quad For any point $(c_{p+1}, \cdots, c_{n-1})$ 
of $\cW$, perform the corresponding change of variables.
Let $pr_2^{-1}(V)\subset \RR^{n-1}$ and, for each $i=1,\cdots p$, 
$(A\cap H_i)_V=(A\cap H_i)
\cap pr^{-1}(V)\subset \RR^{n}$, so that one has 
a commutative diagram 
$$\begin{array}{ccc}
 (A\cap H_i)_V   &\injto  &\RR^n \\
\mapd{pr_1} &     &\mapdr{pr_1} \\
pr_2^{-1}(V)   &\injto &\RR^{n-1}  \\
\mapd{pr_2} &     &\mapdr{pr_2} \\
 V&\injto  &\phantom{\,.}\RR^p\,.
 \end{array}
$$  
Then the fibers of the induced map $pr_1: (A\cap H_i)_V\to pr_2^{-1}(V)$
are finite sets. 
}\smallskip 

\begin{proof}
For each $i$,  $A\cap H_i $ is strictly allowable in $H_i$.
We apply the proposition to $A\cap H_i$, so there exist 
sets $V_i$ and $\cW_i$ with which the conclusion holds for $A\cap H_i$. 
We have only to take $V=\cap V_i$ and $\cW=\cap \cW_i$. 
\end{proof}

We recall here the notion of slicing from [BCR, section 2.3]. 
Let $X=(X_1, \cdots, X_n)$, and let $Y$ be another variable. 
Given a set of polynomials $f_1(X, Y), \cdots, f_s (X, Y)$ in $(X, Y)$, 
[BCR, Theorem 2.3.1] asserts that 
there exists a partition of $\RR^n$ into a finite number of semi-algebraic sets 
$A_1, \cdots, A_m$ and, for each $i=1, \cdots, m$, 
there is a finite number of 
continuous semi-algebraic functions on $A_i$, 
$\xi_{i, 1}<\cdots < \xi_{i, \ell_i}$ such that:
\smallskip 

(i) For each $x\in A_i$, the set $$\{ \xi_{i, 1}(x), \cdots, \xi_{i, \ell_{\ell_i}}(x) \}$$
coincides with the set of roots of those polynomials $f_1(x, Y), \cdots, f_s(x, Y)$
that are not zero polynomials. 

(ii) For each $x\in A_i$, the signs of $f_k(x, y) $, $k=1, \cdots, s$, depend only on the signs of $y-\xi_{i, j}(s)$, $j=1, \cdots, \ell_i$. 

More precisely, this means the following. 
For each $i$ and $j=1, \cdots, \ell_i$, let
$$G(\xi_{i, j})=\{(x, y)\in \RR^{n+1}\mid
x\in A_i\quad\text{and}\quad y=\xi_{i, j}(x)\}$$
be the graph of $\xi_{i, j}$ and, for $j=0, \cdots, \ell_i$,
$$(\xi_{i, j}, \xi_{i, j+1})=\{(x, y)\in \RR^{n+1}\mid
x\in A_i\quad\text{and}\quad \xi_{i, j}(x)<y<  \xi_{i, j+1}(x)\}$$
(by convention, $\xi_{i, 0}=-\infty$ and $\xi_{i, \ell_i+1}=+\infty$). 
Then the sign of $f_k(x, y)$ is constant on each  $G(\xi_{i, j})$
and $(\xi_{i, j}, \xi_{i, j+1})$.
(Here, the sign of a real number $a$ is defined by $\operatorname{sign}(a)=+, 0, -$ according as 
$a>0$, $=0$ or $<0$, respectively. )
\smallskip 

The partition together with a set of functions, 
$(A_i; \xi_{ij})$,
 is called a {\it slicing} of the set of functions $f_1, \cdots, f_s$ with respect to the variable
 $Y$. 

Now assume given a closed semi-algebraic set $A$ of $\RR^{n}=\RR^p\times \RR^{n-p}$
with $p<n$. 
For simplicity, write $r=(r_1, \cdots, r_p)$ and 
$x'=(  x_{p+1}, \cdots, x_{n-1})$. 
Let $(A_i; \xi_{i, j})$ be a slicing, with respect to 
$x_n$,  of a set of functions defining $A$.
Then $A$ is the union of some of the graphs of the functions $\xi_{i, j}$ and 
some of the sets of the form 
$$[\xi_{i, j}, \xi_{i, j+1}]:=\{(r, x', x_n) \in 
\RR^n \mid 
(r, x') \in A_i, \quad \xi_{i, j} (r, x')\le x_n\le \xi_{i, j+1}(r, x')  \}\,.
\eqno{}$$ 
(If $j=0$ interpret the condition as $-\infty <x_n\le \xi_{i, 1}(r, x')$, and similarly for 
$j=\ell_i$.)
\bigskip

Let $A$ be a closed
 semi-algebraic subset of 
$\RR^{p}\times \RR^{n-p}$ with $p<n$. 
We say that  $A$ satisfies {\bf Condition (F)} if
the following holds:
\smallskip 

(F)\quad  For a small $\rho_0>0$, let $V_0=\{(r_1, \cdots, r_p)|\, |r_i|<\rho_0\}$ be an open neighborhood of the origin in $\RR^p$. 
Then for each $i=1, \cdots, p$,  the restriction of the projection $pr_1$, 
$$pr_1: (A\cap H_i)_{V_0} \to pr_2^{-1}(V_0) $$
(see the paragraph preceding (\propfiniteness) for notation) has finite fibers. 
\bigskip 

\noindent (\PropVariation) {\bf Proposition.}\quad{\it
Let $A$ be a (not necessarily closed)
 semi-algebraic subset of 
$\RR^{p}\times \RR^{n-p}$ with $p<n$, satisfying the two conditions:
\smallskip

(i) The closure $\bar{A}$ of $A$ is compact  satisfies Condition (F).

(ii) There are a semi-algebraic set $B$ of $\RR^{n-1}=\RR^p\times \RR^{n-p-1}$ and
continuous semi-algebraic functions $\xi$, 
$\xi'$ defined on $B$, such that $\xi<\xi'$ at each point of $B$ and 
$$A=[\xi, \xi']=\{(r, x', x_n)\in \RR^n|\, 
(r, x')\in B, \quad \xi(r, x')\le x_n\le \xi'(r, x')\}\,.$$
\smallskip 

\noindent We then have: 

(1) For any point $P\in \bar{B}\cap (H_1\cup\cdots\cup H_p)\cap pr_2^{-1}(V_0)$, one has 
$$\lim_{Q\to P} (\xi'-\xi)(Q)=0$$
where the left hand side means the limit as $Q\in B$ tends to $P$. 

(2)  For  $\rho$ with $0<\rho<\rho_0$, let 
$V=\{ |r_i|<\rho\}$
be a smaller neighborhood of the origin.
For each $i=1, \cdots, p$, and for any $\epsilon>0$, there is $\delta$ with $0<\delta<\rho_0$  such that 
$$Q\in B\cap pr_2^{-1}(\bar{V})\quad\text{and}\quad
0< |r_i(Q)| <\delta  \Longrightarrow 
 (\xi'-\xi)(Q) <\epsilon\,.$$
}\smallskip 

\begin{proof}
(1) {\it Step 1.}\quad We show that, for any continuous semi-algebraic function 
$f: [0, 1]\to \RR^{n-1}$ such that $f(0)=P$ and  $f(\, (0, 1])\subset B$, 
we have 
$$\xi\scirc f(t) -\xi'\scirc f(t)\to 0 \quad (t\to 0)\,.$$
For the proof, first note that the compactness of $\bar{A}$
implies that
the functions $\xi$ and $\xi'$ are bounded.  
So the compositions $\xi_1:=\xi\scirc f$ and $\xi'_1:=\xi'\scirc f$
are bounded semi-algebraic functions on $(0, 1]$.
By restricting to an interval $(0, \delta]$ for $\delta>0$ small, and then 
rescaling $t$, one may assume that $\xi_1$ and $\xi'_1$ are continuous on 
$(0, 1]$. 
By [BCR, Proposition 2.5.3], these functions extend to continuous semi-algebraic 
functions on $[0, 1]$, which we still denote by the same $\xi_1$ and $\xi'_1$. 
The set 
$$[\xi_1, \xi'_1]:=\{(f(t), x_n)|\,  0\le t\le 1, \quad \xi_1(t) \le x_n\le \xi'_1(t)\}$$
is the closure of the set
$$\{(f(t), x_n)|\,  0< t\le 1, \quad \xi_1(t) \le x_n\le \xi'_1(t)\}\,,$$
and the latter is a subset of $A$;
thus  $[\xi_1, \xi'_1]$ is contained in $\bar{A}$. 
Since the fiber of $P$ of the projection from $[\xi_1, \xi'_1]$ to $\RR^{n-1}$ is 
the interval $[\xi_1(0), \xi'_1(0)]$, and since $\bar{A}$ 
satisfies (F), one must have 
$\xi_1(0)=\xi'_1(0)$.

{\it Step 2.}\quad We show the assertion of (1). 
If it were false, there exists $\eps>0$ such that the semi-algebraic set 
$$B_\eps= \{Q\in B| (\xi' -\xi)(Q)\ge \eps\}$$
contains in its closure the point $P$. 
By the curve selection lemma ([BCR, Theorem 2.5.5]) there is a continuous 
semi-algebraic map $f: [0, 1]\to \RR^{n-1}$ such that $f(0)=P$
and $f(\, (0, 1])\subset B_\eps$.  One then has
$\xi'\scirc f(t) -\xi\scirc f(t)\ge \eps$ for any $t>0$. 
This contradicts the assertion proven in Step 1. 

(2)  Assume that the assertion was false. There is $\epsilon>0$ and a sequence of points $Q_n\in B\cap pr_2^{-1}(\bar{V})$ with 
$r_i (Q_n)\to 0$ and $(\xi'-\xi)(Q) \ge \epsilon$. 
By compactness of $\bar{B}$, there is a convergent subsequence $Q_{n'}$, which 
converges to a point $P\in \bar{B}$. Then one necessarily has 
$P\in \bar{B}\cap pr_2^{-1}(V_0)\cap \{r_i=0\}$. 
But by  (1), there is $\delta>0$ such that for any 
$Q\in B$ with $|Q-P|<\delta$, one has  $(\xi'-\xi)(Q) < \epsilon$.
This is a contradiction. 
\end{proof}

Let $A$ be a semi-algebraic subset of 
$\RR^{p}\times \RR^{n-p}$ with $p<n$ such that $\bar{A}$ is compact. 
Define a function on $(r_1, \cdots, r_p)\in \RR^p$ by
$$v(r_1, \cdots, r_p)=v(A|x_n; (r_1, \cdots, r_p)\,)=
\sup_{ T \mapsto (r_1, \cdots, r_p)}
 vol(A\cap pr_1^{-1}(T); dx_n)$$
where $T\in \RR^{n-1}$ varies over the points mapping to $(r_1, \cdots, r_p)$ 
by the projection $pr_2: \RR^{n-1} \to \RR^p$, and
$pr_1:\RR^n \to \RR^{n-1}$ is the first projection. 
Since $\bar{A}$ is compact, $v(A|x_n; (r_1, \cdots, r_p)\,)$ takes finite
values, so $v$ is a bounded semi-algebraic function in $(r_1, \cdots, r_p)$.
We will give an estimate for this function, assuming Condition (F) on $\bar{A}$. 
\bigskip

\noindent (\PropLtypeineq) {\bf Proposition.}\quad{\it
Let $A$ be a semi-algebraic subset of 
$\RR^{p}\times \RR^{n-p}$ with $p<n$, such that $\bar{A}$ is compact 
and satisfies Condition (F).
Then there exist $0<\rho<Min(\rho_0, 1)$ and 
$C, \al >0$ such that for any 
$(r_1, \cdots, r_p)$ with $0<|r_i|<\rho$, one has
$$v(A|x_n; (r_1, \cdots, r_p)\,)\le C |r_1|^\al |r_2|^\al \cdots |r_p|^\al\,.$$}

\begin{proof}
 Taking a slicing $(A_i; \xi_{i,j})$ of $A$, we can assume that $A$ is either the graph
$G(\xi_{i, j})$ of a function $\xi_{i, j}$ or a set
of the form $[\xi_{i, j}, \xi_{i, j+1}]$ for $1\le j\le \ell_{i}-1$. 
In the former case, $v(r_1, \cdots, r_p)=0$, so the assertion is obvious.
We will thus assume $A$ is of the latter form. 

Consider the function  
on $|r_1|<\rho$ defined by
$$v_1(r_1)=\sup_{(r_2, \cdots, r_p)\in \RR^{p-1}} 
v(r_1, r_2, \cdots, r_p) \,,$$
where $(r_2, \cdots, r_p)$ varies over points with $|r_i|<\rho$. 
This function is semi-algebraic and bounded.  
If $\rho$ is taken small enough with $0<\rho<Min(\rho_0, 1)$, the function $v_1(r_1)$ is 
continuous and semi-algebraic on $0<|r_1|<\rho$. 
Similarly the function $v_i(r_i)$, $i=1, \cdots, p$, is defined taking the supremum of the values
of $v$ with $r_i$ fixed, and it is continuous semi-algebraic on  $0<|r_i|<\rho$. 

We claim that $\lim_{r_1\to 0} v_1(r_1)=0$. 
Indeed by (2) of (\PropVariation), for any $\eps>0$, there exists $\delta>0$ such that for 
any $Q\in\RR^{n-1}$ with $|r_i(Q)|<\rho$ and $0<|r_1(Q)|<\delta$, one has 
$vol(A\cap pr_1^{-1}(Q); dx_n)<\eps$.  Thus if $0<|r_1|<\delta$, then 
$v_1(r_1)\le \eps$. 

Setting $v_1(0)=0$, the function $v_1(r_1)$ is continuous and semi-algebraic  on $|r_1|<\rho$.  
Applying \L{}ojasiewicz's inequality ([BCR, Corollary 2.6.7]), 
there exist positive numbers
$C, \al$ such that 
$$v_1(r_1)\le C |r_1|^\al\,$$
for $|r_1|<\rho$. 
By the same reasoning, one obtains inequalities (with the same $C, \al$, 
by taking $C$ larger and $\al$ smaller if necessary)
$v_i(r_i)\le C |r_i|^\al$ for $0<|r_i|<\rho$, for each $i=1, \cdots, p$. 
Combined with the evident inequalities $v(r_1, \cdots, r_p)\le v_i(r_i)$, one obtains 
$$v(r_1, \cdots, r_p)\le (v_1(r_1)\cdots v_p(r_p)\,)^{1/p}
\le C |r_1|^{\al/p} \cdots| r_p|^{\al/p}\,.$$
\end{proof}

Before stating the next proposition, we introduce some notation. 
Let $\mu: \ti{\RR}^n \to \RR^n$ be a succession of permissible blow-ups.
A face of $\RR^n=\RR^r\times \RR^{n-p}$ is of the form $(\mbox{face of}\times \RR^p)\times
\RR^{n-p}$; similarly on a permissible blow-up of $\RR^n$.
Thus 
one has a succession of permissible blow-ups $\bar{\mu}: \ti{\RR}^p \to \RR^p$
such that $\ti{\RR}^n=\ti{\RR}^p\times \RR^{n-p}$ and 
$\mu$ is the product of $\bar{\mu}$ 
and the identity map on $\RR^{n-p}$.
Note that in $\RR^p$ or its blow-up, minimal faces are points, which we may call {\it vertices}.
Letting  $ \mu'=\bar{\mu}\times 1 :\ti{\RR}^{n-1}:=\ti{\RR}^p\times \RR^{n-p-1} \to \RR^p\times \RR^{n-p-1}=\RR^{n-1}$,
 we obtain the following commutative diagram 
$$\begin{array}{crc}
\RR^n&\mapl{\mu}&\ti{\RR}^n \\
\mapd{pr_1}&&\mapdr{pr_1} \\
\RR^{n-1}&\mapl{\mu'}&\ti{\RR}^{n-1}\\
\mapd{pr_2}& &\mapdr{pr_2} \\
\RR^p&\mapl{\bar{\mu}}&\ti{\RR}^p
\end{array}$$
where $pr_1$, $pr_2$ on the left are the projection maps 
as introduced before, and the maps $pr_1$, 
$pr_2$ on the right are similarly defined projection maps. 
Also let $pr=pr_2pr_1:\ti{\RR}^n \to \ti{\RR}^p$. 
\bigskip 

\noindent (\PropBlupfinite) 
{\bf Proposition.}\quad{\it  
Let $A$ be an allowable 
compact semi-algebraic subset of 
$\RR^{p}\times \RR^{n-p}$ with $p<n$. 
Then there exist a succession of permissible
blow-ups
$\mu: \tilde{\RR}^n\to \RR^n$ 
and  a point $(c_{p+1}, \cdots, c_{n-1})\in \RR^{n-p-1}$
such that, after the corresponding change of variables,
one has the following property:
\smallskip 

There is an open set $V$ of $\ti{\RR}^p$ containing all the vertices
such that, for each $i=1, \cdots, p$, 
the induced map 
$$pr_1: (\mu^{-1}(A)\cap H_i)_{V}\to pr_2^{-1}(V)$$
has finite fibers.
}\smallskip 

\begin{proof} By Theorem (\thmreduction) one can take a succession of permissible blow-ups
$\mu: \ti{\RR}^n\to \RR^n$ such that $\mu^{-1}(A)$ is almost strictly allowable. 
Then apply Corollary (\PropVariation.1) to each chart of $\ti{\RR}^n$ (which corresponds to a vertex of
$\ti{\RR}^p$). The intersection of the open sets $\cW$ for the vertices is non-empty, 
so there exists $(c_{p+1}, \cdots, c_{n-1})$ as asserted.
\end{proof}

\noindent (\Thmconvergence) {\bf Theorem.}\quad{\it
Let $A$ be a compact semi-algebraic subset 
of $\RR^p\times \RR^{n-p}$ (with $0\le p\le n$) which is 
allowable, and let $\omega$ be
a logarithmic $n$-form.  Then the integral 
$$\int_A \omega$$
absolutely converges. 
}\bigskip  

By compactness of $A$, it is enough to show the following local version of the theorem: 
\bigskip 

\noindent (\Thmconvergence.1) {\bf Theorem}\quad{\it (Under the same assumption as for (\Thmconvergence)\,.)
For any point $P$ of $A$, there exists a neighborhood $W$ of
$P$ in $\RR^n$ such that the integral
$$\int_{A\cap W} \omega\,$$
 is absolutely convergent. 
}\smallskip 

\begin{proof}
We prove this by induction on the cardinality $p$ of the set $\{i \,|\,r_i(P)=0\}$. 
If it is zero, the form is smooth in an neighborhood of $P$, so the assertion is obvious.
If $p$ is positive, by translation in the $x$-coordinates 
we may assume that $P$ is the origin in $\RR^p\times \RR^{n-p}$. 
Note that $p<n$ by allowability. 

Take a succession of permissible blow-ups $\mu$ as in
 Proposition (\PropBlupfinite).
With this, we claim:
\smallskip 

(\Thmconvergence.2) {\bf Claim.}\quad {\it The integral $\int_{\mu^{-1}(A)} \mu^*\omega
$
is absolutely convergent in a neighborhood of
each point $Q$ of $\mu^{-1}(A)$. 
}\smallskip 

Indeed, if $Q$ does not lie over a vertex, the number of codimension one faces of $\ti{\RR}^n$
containing $Q$ is less than $p$, so by induction hypothesis for (\Thmconvergence.1) we have absolute convergence.
Assuming that $Q$ lies over a vertex of $\ti{\RR}^p$, 
let $(r'_1,\cdots, r'_p)$ be the coordinates for
the chart $\RR^p$ at $pr(Q)$. 
After a change of coordinates, 
(\PropLtypeineq) is satisfied on $V=\{|r'_i|<\rho\}$. 
Let $W=pr^{-1}(V)$, and rewrite $(r_1, \cdots, r_p)$
for $(r'_1,\cdots, r'_p)$.
We may assume that 
$\omega={dr_1}/{r_1}\wedge\cdots\wedge {dr_p}/{r_p}\wedge dx_{p+1}\wedge\cdots\wedge dx_n$
on $W$.
By compactness $A$ is contained in the set $\{ a\le x_i\le b |\, i=p+1, \cdots, n-1\}$, for some $a<b$; we 
set $M=(b-a)^{n-p-1}$. 
Applying Fubini's theorem to $\mu^{-1}(A)\cap W\subset \RR^n =\RR^{n-1}\times \RR$, 
then to $V\times [a, b]^{n-p-1}\subset \RR^{n-1}$, and taking (\PropLtypeineq) into account, 
one has
\begin{eqnarray*}
&&vol(\mu^{-1}(A)\cap W; \frac{dr_1}{r_1}\wedge\cdots\wedge \frac{dr_p}{r_p}\wedge dx_{p+1}\wedge\cdots\wedge dx_n) \\
&=& \int_{\mu^{-1}(A)\cap W} \frac{dr_1}{|r_1|}\cdots \frac{dr_p}{|r_p|}\,dx_{p+1}\cdots dx_{n-1}dx_n \\
&= &\int_{V\times [a, b]^{n-p-1}}vol (\mu^{-1}(A)\cap pr_1^{-1}(T); dx_n)\,\frac{dr_1}{|r_1|}\cdots \frac{dr_p}{|r_p|}\,dx_{p+1}\cdots dx_{n-1} \\
&\le&M\int_V  v(A|x_n; (r_1, \cdots, r_p)\,)
\,\frac{dr_1}{|r_1|}\cdots \frac{dr_p}{|r_p|} \\
&\le&M\int_{0<|r_i|<\rho} C |r_1|^\al |r_2|^\al \cdots |r_p|^\al 
\,\frac{dr_1}{|r_1|}\cdots \frac{dr_p}{|r_p|} < +\infty\,,
\end{eqnarray*} 
since $\int_{0<r<\rho} r^\al \frac{dr}{r} <+\infty$. 

We have only to show that (\Thmconvergence.2) implies (\Thmconvergence.1) for $A$ at 
$P=0$.  For each $Q$ of $\mu^{-1}(A)$, take a neighborhood $W_Q$ such that 
$vol(\mu^{-1}(A)\cap W_Q, \mu^*\omega)$ is absolutely convergent. 
Since $\mu^{-1}(A)$ is compact, there is a finite number of them, $W_{Q_k}$, 
that covers $\mu^{-1}(A)$. 
Then 
$$vol(A, \omega)=vol(\mu^{-1}(A), \mu^*\omega)\le 
\sum vol(\mu^{-1}(A)\cap 
W_{Q_k},  \mu^*\omega)<+\infty.$$
\end{proof}

In the rest of this section, we take a monomial $u=r_1^{m_1}\cdots r_p^{m_p}\neq 1$,
and for $A$ an allowable semi-algebraic subset of 
$\RR^p\times \RR^{n-p}$, 
consider the set $A\cap \{u=t\}$, where $t\in \RR-\{0\}$ with $|t|$ small.
Since $u$ is not a constant function, this is of dimension $\le n-1$.
 For a logarithmic $(n-1)$-form $\omega$, 
we will study the asymptotic behavior of the integral of $\omega$ over 
$A\cap \{u=t\}$. 
\bigskip 

\noindent (\Propnormallink) {\bf Proposition.}\quad{\it 
Let $A$ be a semi-algebraic set of $\RR^n$ such that its closure $\bar{A}$ is compact and
satisfies Condition (F), and assume 
$\rho>0$ is chosen so that the statement of (\PropLtypeineq)
holds.  
 Then for any logarithmic $(n-1)$-form $\omega$, and any monomial $u=r_1^{m_1}\cdots r_p^{m_p}\neq 1$, 
there exist $C', \al'>0$ such that  
$$vol \bigl(\,(A\cap W_\rho
 \cap \{u=t\}) ;\, \omega\bigr) 
\le C'\,|t|^{\al'}\,.$$
(We set $W_\rho:=\{0<|r_i|<\rho\, \, (i=1, \cdots, p)\}\subset \RR^p\times \RR^{n-p}$.)
}\smallskip 

\begin{proof}
For the proof we consider a variant of the function 
$v(A|x_n; (r_1, \cdots, r_p)\,)$ introduced before (\PropLtypeineq). 
Let $pr_1: \RR^{n}\to \RR^{n-1}$ be the projection as before, and 
$pr_3:\RR^{n-1}\to \RR^{p-1}$ be the projection to 
the coordinates $(r_2, \cdots, r_p)$. 
For an algebraic set $B\subset \RR^n$, we set  
$$v(B \,|x_n; (r_2, \cdots, r_p)\,)=
\sup_{ T \mapsto (r_2, \cdots, r_p)}
vol(B\cap  pr_1^{-1}(T)\,; dx_n)$$
where $T\in \RR^{n-1}$ varies over the points mapping to $(r_2, \cdots, r_p)$ 
by the projection $pr_3$. 
If $B=A\cap W_\rho\cap \{u=t\}$ we have:
\bigskip 

\noindent (\Propnormallink.1) {\bf Lemma.}\quad{\it Under the same assumption as for (\Propnormallink),
there exist $C, \al'>0$ such that
$$v\bigl(\,(A\cap W_\rho\cap \{u=t\}) \,|x_n; \,(r_2, \cdots, r_p)\,\bigr)
\le C|t|^{\al'}
|r_2|^{\al'}\cdots |r_p|^{\al'}\,$$
for $(r_2, \cdots, r_p)$ with $0<|r_i|< \rho$, $2\le i\le p$, and $t\in \RR-\{0\}$ with $|t|$ small.
}
\begin{proof}
Note that the set 
$$A\cap W_\rho\cap\{u=t\}\cap  pr_1^{-1}(T)$$
 equals $A\cap pr_1^{-1}(T)$ if $pr_2(T)=(r_1, \cdots, r_p)$ 
 satisfies $u=t$ and $0<|r_i|<\rho$; otherwise it is empty
 (recall that $pr_2: \RR^{n-1}\to \RR^p$ is projection by $(r_1, \cdots, r_p)$). 
Therefore
\begin{eqnarray*}
&&v(A\cap W_\rho\cap \{u=t\}|x_n; (r_2, \cdots, r_p)\,)  \\
&=&\sup_{r_1} v(A|x_n; (r_1, \cdots, r_p)\,)\quad
(\text{supremum taken over $r_1$ such that $u=t$ and $0<|r_1|<\rho$})\\
&\le& \sup_{r_1} C|r_1|^\al |r_2|^{\al}\cdots |r_p|^\al\,.
\quad(\text{by (\PropLtypeineq)})
\end{eqnarray*}

It is thus enough to show that there exists $\al'>0$ such that
$$|r_1|^\al |r_2|^{\al}\cdots |r_p|^\al\le |u|^{\al'} |r_2|^{\al'}\cdots |r_p|^{\al'}\,.
$$
Since the right hand side equals 
$|r_1|^{m_1\al'}|r_2|^{(m_2+1)\al'}
\cdots |r_p|^{(m_p+1)\al'} $ and 
$|r_i|<1$, it suffices to take 
$\al'>0$ small enough so that
$$m_1\al'\le \al, \quad (m_2+1)\al'\le \al,\, \,\cdots, \quad (m_p+1)\al'\le \al\,.$$
\end{proof}

We return to the proof of (\Propnormallink). 
On $\{u=t\}$, there is a relation
$\sum_{1\le i\le  p} m_i {dr_i}/{r_i}=0$,
hence $dr_1/r_1\wedge\cdots\wedge dr_p/r_p=0$.
So the restriction of the form $\omega$ 
can be written 
$$\omega|_{\{u=t\}}
=\sum_{i=1}^{n} a_i \frac{dr_1}{r_1}\wedge\cdots\wedge 
\frac{dr_{i-1}}{r_{i-1}}\wedge
\frac{dr_{i+1}}{r_{i+1}}\wedge\cdots\wedge
\frac{dr_p}{r_p}\wedge dx_{p+1}\wedge\cdots\wedge dx_n$$
with some polynomial functions $a_i$ on $\{u=t\}$.
Thus it is enough to prove the assertion for the case
 $\omega=\frac{dr_2}{r_2}\wedge\cdots\wedge \frac{dr_p}{r_p}\wedge
dx_{p+1}\wedge\cdots\wedge dx_n$. 
Assume that $A$ is contained in the set $\{ a\le x_i\le b |\, i=p+1, \cdots, n-1\}$ and 
set $M=(b-a)^{n-p-1}$. 
Then one has 
\begin{eqnarray*}
&&vol(A\cap W_\rho \cap \{u=t\}; \frac{dr_2}{r_2}\wedge\cdots\wedge \frac{dr_p}{r_p}\wedge
dx_{p+1}\wedge\cdots\wedge dx_n)\\
&\le & M\int_{\{0<|r_i|\le  \rho\}} C |t|^{\al'} |r_2|^{\al'}\cdots |r_p|^{\al'} 
\frac{dr_2}{|r_2|}\wedge\cdots\wedge \frac{dr_p}{|r_p|} \\
&=& C'|t|^{\al'}\,.
\end{eqnarray*}

\end{proof}

\noindent (\Thmnormallink) {\bf Theorem.}\quad{\it Let $A$ be a compact semi-algebraic subset 
of $\RR^p\times \RR^{n-p}$ (with $0\le p\le n$) which is allowable. 
Let $u=r_1^{m_1}\cdots r_p^{m_p}\neq 1$ be a monomial, and $\omega$ be a logarithmic $(n-1)$-form.
Then there exist $C, \al>0$ such that one has, for $t\in \RR$ with $|t|\neq 0$ small,
$$vol(A \cap \{u=t\}; \omega)\le
C\,|t|^\al\,.$$
}\bigskip

As for (\Thmconvergence), we can reduce it to 
a local statement:
\bigskip 

\noindent (\Thmnormallink.1) {\bf Theorem.}\quad{\it
(Under the same assumption.) 
For any point $P$ of $A$, there exists an open neighborhood $W$ of 
$P$ in $\RR^n$ such that for some $C, \al>0$, one has an inequality 
$$vol(A\cap W\cap \{u=t\}; \omega)\le
C|t|^\al\,.$$}
\smallskip 

\begin{proof}
The proof proceeds by induction on $p=\sharp \{i|\, r_i(P)=0\}$. 
If $p=0$, the set $A\cap \{u=t\}$ does not contain $P$, so the assertion obviously holds.
One may thus assume that $P$ is the origin in $\RR^p\times \RR^{n-p}$ with $p\ge 1$; we have $p<n$ by allowability.

Take a succession of permissible blow-ups $\mu$ as in
 (\PropBlupfinite).
\smallskip 

(\Thmnormallink.2) {\bf Claim.}\quad {\it 
For any point $Q$ of $\mu^{-1}(A)$, there exist an open neighborhood $W_Q$ of $Q$
and $C, \al>$ such that 
$$vol(A\cap W_Q\cap \{u=t\}; \mu^* \omega)\le
C|t|^\al\,.$$
}\smallskip 

Indeed, if $Q$ does not lie over a vertex, the assertion holds 
by induction hypothesis for (\Thmnormallink.1).
If $Q$ lies over a vertex of $\ti{\RR}^p$, 
let $(r'_1,\cdots, r'_p)$ be the coordinates for
the chart at $pr(Q)$. 
After a change of coordinates, 
(\Propnormallink) is satisfied on $V=\{|r'_i|<\rho\}$. 
Then the assertion follows from (\Propnormallink).

Finally, we derive (\Thmnormallink.1) from (\Thmnormallink.2)
as in the proof of (\Thmconvergence).
\end{proof}

For reference in \Part II, we record a special case. 
\bigskip
 
\noindent (\Thmnormallink.3) {\bf Corollary.}\quad{ In particular, taking $u=r_1$, we have
$$vol(A\cap \{r_1=t\};\frac{dr_2}{r_2}\wedge\cdots\wedge \frac{dr_p}{r_p}\wedge
dx_{p+1}\wedge\cdots\wedge dx_n)\le C\,|t|^\al\,.$$
}
\bigskip

\newpage

\newcommand{\BHZ}{{\BH}^0}
\newcommand{\HZ}{H^0}

{\bf \S \sectadmiss. Integrals on admissible semi-algebraic sets}
\bigskip

Let $\CC^n$ be affine $n$-space over $\CC$, with coordinates $(z_1, \cdots, z_n)$. 
Let $\BH_i=\{z_i=0\}$, $i=1, \cdots, n$,  be  the coordinate hyperplanes.
An intersection of coordinate hyperplanes will be called a face; it is of the form
$\BH_I=\cap_{i\in I} \BH_i$ 
for a subset $I$ of $\{1, \cdots, n\}$.   
The union of the divisors $\{z_i=1\}$ for $i=1, \cdots, n$ will be denoted by $D$. 

Let $(\BP^1)^n$ be the $n$-fold product of projective line $\BP^1$, with
 $(z_1, \cdots, z_n)$ affine coordinates. 
For $i=1, \cdots, n$ and $\al\in \{0, \infty\}$, let
$\BH_i^\al=\{z_i=\al\}$ be a prime divisor; by definition a face of $(\BP^1)^n$
is the intersection of some of these divisors.
Let  $D$ be the union of the divisors $\{z_i=1\}$, and
$\sq^n=(\BP^1)^n-D$ be the complement of $D$.

There are canonical isomorphisms $\BP^1-\{0\} \cong \CC$, $(1: z) \leftrightarrow z$, and 
$\BP^1-\{\infty\} \cong \CC$, $(z: 1) \leftrightarrow z$.
Thus given a sequence $(\al_1, \cdots , \al_n)$ with $\al_i\in \{0, \infty\}$, 
there is an isomorphism $\prod (\BP^1-\{\al_i\}) \cong \CC^n$; we refer to these as the charts of $(\BP^1)^n$.  On each chart, the normal crossing divisor $\sum \BH_i^\al $ restricts 
to the divisor $\sum \BH_i$, and the divisor $D$ restricts to $D$.
\bigskip 

\noindent (\Dadmissible) {\bf Definition.}\quad 
Let $m$ be an integer with $n\le m\le 2n$.
A closed semi-algebraic set $A$ of $\CC^n$
 is said to be {\it $m$-admissible}
(with respect to $\{\BH_i\}$) if for any set $\BH_I$ one has
$$\dim (A\cap \BH_I\cap (\CC^n-D)\, )\le m-2|I|\,.$$
(In particular, $\dim A\le m$).

Similarly, a closed semi-algebraic set $A$ of $(\BP^1)^n$
 is said to be $m$-admissible (with respect to $\{\BH_i^\al\}$) if 
for any face $F$ one has
$$\dim (A\cap F \cap \sq^n)\le m- 2\operatorname{codim}_\CC (F)\,$$
\bigskip 

If $A$ is an $m$-admissible closed semi-algebraic set of $\CC^n$ (resp. $(\BP^1)^n$), 
then any closed semi-algebraic subset of $A$ is also $m$-admissible. 
If $A$ is an $m$-admissible closed semi-algebraic set of $\CC^n$,
then $A\cap \BH_i$ is $(m-2)$-admissible in $\BH_i=\CC^{n-1}$. 
If $A$ is  $m$-admissible in $(\BP^1)^n$, then for each chart $U\cong \CC^n$,
$A\cap U$ is $m$-admissible in $\CC^n$.
\bigskip

\noindent (\Spolar) {\bf Change of variables via the polar coordinates.}\quad
We will take an $m$-admissible compact semi-algebraic set $A$ of $\CC^n$ (or $(\BP^1)^n$)
and a logarithmic $(n, m-n)$-form $\omega$ on $\CC^n$, and study the integral
$\int_A \omega$. 
First assume $A$ is a subset of $\CC^n$.  

We will introduce a change of variables as follows.
In case $n=1$, let 
$$S_i=\{z=re^{i\theta}\, |\, 0\le r<+\infty, \quad  (i-2)\pi/4 \le \theta \le i\pi/4\,\}\,,$$
for $i=1,2,3, 4$, be sectors which cover the complex plane.
For simplicity we put 
 $\CC_+=S_1$. 

On the sector $S_1=\CC_+$, let 
$$r=|z|, \quad \tau=y/x \quad\text{ for }\quad
 z=x+i y \,.$$ 
Note that if $(r, \theta)$ are the polar coordinates for 
$z$, then $\tau=\tan \theta$. 
The map 
$\pi: \RR_{\ge 0}\times [-1, 1] \to \CC
$
which sends $(r, \tau)$ to $x+iy$ with
$$x= \frac{r}{\sqrt{\tau^2+1}}, \qquad y=
\frac{r \tau }{\sqrt{\tau^2+1}}$$
induces a continuous semi-algebraic map 
$$\pi: \RR_{\ge 0}\times [-1, 1] \to 
\CC_+ \,, $$
which is proper surjective and isomorphic outside the origin of
the sector.
Note that the map by the polar coordinates
$(r, \theta)\mapsto (x, y)$ is {\it not} semi-algebraic, 
but by using $\tau$ in place of $\theta$ we obtain the semi-algebraic map $\pi$. 

 On $\CC_+ -\{0\}$, we have identities  
$$\frac{dz}{z}=\frac{dr}{r}+ i\frac{d\tau}{\tau^2+1}\,$$
and 
$$\frac{dz}{z}\wedge d\bar{z}=\frac{-2(\tau+i)}{(\tau^2+1)^{3/2}}\,dr\wedge d\tau\,.$$
For $S_2$, take  $\tau=-x/y$ and consider the map 
$\pi: \RR_{\ge 0}\times [-1, 1] \to S_2$
that sends $ (r, \tau)$ to  $-{r \tau }/{\sqrt{\tau^2+1}} +{i \tau }/{\sqrt{\tau^2+1}} $.
Then the same formula for ${dz}/{z}\wedge d\bar{z}$ holds up to a factor of $-i$. 
Similarly for $S_3$ and $S_4$, with appropriate choices of $\tau$ and $\pi$. 

In case $n\ge 2$,  $\CC^n$ is covered by the products of sectors
$$S_{\al_1}\times \cdots \times S_{\al_n}$$
for sequences $(\al_1, \cdots, \al_n)$ taking values in 
$\{1, 2, 3, 4\}$. 
For each such product of sectors, there is a continuous, proper surjective  semi-algebraic
map
$$(\RR_{\ge 0}\times [-1,1])^n \to S_{\al_1}\times \cdots \times S_{\al_n}$$
given as the product of the maps $\pi:\RR_{\ge 0}\times [-1,1]\to S_{\al_i}$ just mentioned.
\bigskip 
 
Let $m$ be an integer with $n\le m\le 2n$. 
By definition, an $(n, m-n)$-form on $\CC^n$ with logarithmic singularities along $\BH_i$ is 
 a linear combination, with coefficients smooth 
functions on $\CC^n$, of the forms
$$(\frac{dz_1}{z_1}\wedge \cdots \wedge \frac{dz_n}{z_n})\wedge
\bigwedge_{k\in R} d\bar{z}_k\eqno{(\Spolar.a)}$$
with $R$ a subset of cardinality $m-n$ of $\{1, \cdots,n\}$. 

 Writing $\ti{\CC}_+=\RR_{\ge 0}\times [-1, 1]$ for short,  let 
$$\pi=\pi_{\al_1, \cdots, \al_n}: \ti{\CC}_+^n \to S_{\al_1}\times \cdots \times S_{\al_n}\eqno{(\Spolar.b)}$$
be the product of the maps $\pi:\ti{\CC}_+ \to S_{\al_i}$, still denoted by the same letter. 
It is a continuous, proper surjective semi-algebraic map, which is
 isomorphic over $\prod (S_{\al_i}-\{0\})$. 
The map $i\scirc \pi:  \ti{\CC}_+^n \to \CC^n$, where $i: S_{\al_1}\times \cdots \times S_{\al_n}\to \CC^n$
is the inclusion, will also be written $\pi$.

The pull-back of the form (\Spolar.a) by $\pi$ is  the sum of 
the forms (up to a product of $\pm i$)
$$(\prod_{j\in Q} \frac{1}{\tau_j^2+1}) \left(\prod_{k\in R}
\frac{-2(\tau_k+i )}{(\tau_k^2+1)^{3/2}}\right)
\cdot\vphi_{P, Q, R}\,$$
with
$$\vphi_{P, Q, R}:=\bigwedge_{i\in P}\frac{dr_i}{r_i}
\wedge \bigwedge_{j\in Q} d\tau_j\wedge\bigwedge_{k\in R}(dr_k \wedge
d\tau_k )\,,\eqno{(\Spolar.c)}$$
where $(P, Q)$ varies over the partitions of $\{1, \cdots, n\}-R$.
The function in front of $\vphi_{P, Q, R}$ in the first formula is smooth and bounded.

Consider the map
$$\ti{\CC}_+^n = \ti{\CC}_+^P\times \ti{\CC}_+^Q\times \ti{\CC}_+^R\to \RR_{\ge 0}^P \times [-1,1]^Q 
\times \ti{\CC}_+^R$$
obtained as the product of the maps
$$\begin{array}{ll}
\ti{\CC}_+^P \to \RR_{\ge 0}^P\,,&(r_i, \tau_i) \mapsto (r_i)\,, \\
\ti{\CC}_+^Q \to [-1, 1]^Q\,,&(r_j, \tau_j)\mapsto (\tau_j)\,,\\
\ti{\CC}_+^R \to \ti{\CC}_+^R\,, &\text{the identity map\,.}
\end{array}
$$
The target of this map is a subset of 
$\RR^P\times \RR^Q\times (\RR\times \RR)^R$
with coordinates $(r_i)_{i\in P}$, $(\tau_j)_{j\in Q}$ and $(r_k, \tau_k)_{k\in R}$.
Composing with the inclusion we obtain a map
$$q=q_{P, Q, R}: \ti{\CC}_+^n \to \RR^P\times \RR^Q\times (\RR^2)^R=\RR^p\times \RR^{m-p}\,.$$
For $i\in P$, let $H_i=\{r_i=0\}$, and for $I\subset P$, $H_I=\cap_{i\in I}H_i$. 
\bigskip 

For $A$ an $m$-admissible compact semi-algebraic set of $\CC^n$, and $\omega$ a
logarithmic $(n, m-n)$-form, there exists an open semi-algebraic set $U\subset A$ which
is a Nash submanifold and  $\dim (A\backslash U)<m$, thus $\int_A|\omega|$ is defined.
We ask if the integral is absolutely convergent.  Taking intersections, we may
assume that $A$ is contained in $S_{\al_1}\times \cdots \times S_{\al_n}$. 
Let $\pi$ be the map (\Spolar.b). 

If $A^0:= A-(\cup \BH_i)$, then by definition $\int_A |\omega|$ is equal to 
 $$\int_{A^0} |\omega|\,.$$
 Since $\pi: \pi^{-1}(A^0) \to A^0$ is an isomorphism of semi-algebraic sets,
 it equals to $\int_{\pi^{-1}(A^0)} |\omega|$.
Decomposing the form $\omega$ as above,
 one is reduced to the absolute convergence of the integral
$$\int_{\pi^{-1}(A^0)} |\vphi_{P, Q, R}|\,.$$

Let $\vphi'_{P, Q, R}$ be the $m$-form on $\RR^P\times \RR^Q\times (\RR^2)^R$
defined by 
the same formula (\Spolar.c), which has logarithmic singularities along 
$\{r_i=0\}$ for $i\in P$.
We have 
$\vphi_{P, Q, R}=q^*\vphi'_{P, Q, R}$.
Applying \APropboundintegral 
to the set $\pi^{-1}(A^0)$ and the projection $q$, 
we further reduce the question to the absolute convergence of 
$$\int_{q\pi^{-1}(A^0)} |\vphi'_{P, Q, R}|\,.
\eqno{(\Spolar.d)}
$$

The last convergence holds true if the closure of $q\pi^{-1}(A^0)$ is 
allowable in $\RR^P\times \RR^Q\times (\RR^2)^R$
with respect to $\{r_i=0\}$, $i\in P$. This will be shown in the next proposition.

Next assume that $A$ is compact $m$-admissible in $(\BP^1)^n$. 
For a chart $\CC^n=\prod (\BP^1-\{\al_i\})^n$,
let  $R_i=\{|z_i|\le 1\}$ if $\al_i=\infty$ and $R_i=\{|z_i|\ge 1\}$ if $\al_i=0$, and
 consider the subset of the chart
$$A'= A\cap (R_1\times \cdots\times R_n)\,.$$
Then $A'$ is compact $m$-admissible in the chart. 
Since the integral $\int_{A'}|\omega|$ absolutely converges
on each chart,  the integral $\int_A |\omega|$ also absolutely converges.
\bigskip

\noindent (\admissallow) {\bf Proposition.}\quad{\it  
Let $A$ be a compact semi-algebraic set of $\CC^n$
which is $m$-admissible with respect to $\{\BH_i\}$.
Let $\pi:\ti{\CC}_+^n \to \CC^n$ be a map as (\Spolar.b), and
$q: \ti{\CC}_+^n \to \RR^m$ be the map associated with a partition $(P, Q, R)$ of 
$\{1, \cdots, n\}$ with $|P|+|Q|+2|R|=m$.
Then the set 
$q \pi^{-1}(A) $
is a compact semi-algebraic set of $\RR^m$, which is allowable
with respect to $H_i=\{r_i=0\}$, $i\in P$. }
\smallskip 

\begin{proof}  We begin by introducing auxiliary maps $\pi_1$ and $q_1$.
For simplicity assume that $\pi=\pi_{1, \cdots, 1}: \ti{\CC}_+^n \to \CC_+^n$. 
For the partition $(P, Q, R)$, 
consider the product map 
$$1\times \pi\times \pi: \CC^P \times \ti{\CC}_+^Q\times  
\ti{\CC}_+^R\to 
\CC^P \times \CC_+^Q\times  
\CC_+^R\,, $$
and let $\pi_1:\CC^P \times \ti{\CC}_+^Q\times  
\ti{\CC}_+^R\to \CC^n$ be 
the composition with the inclusion in $\CC^n$. 
The map $\pi: \ti{\CC}_+^n\to \CC^n$ factors as $\pi_1\scirc pr$,
where $pr: \ti{\CC}_+^n\to \CC^P \times \ti{\CC}_+^Q\times  
\ti{\CC}_+^R$ is projection. 
Also consider the map
$$q_1: {\CC}_+^P\times \ti{\CC}_+^Q\times \ti{\CC}_+^R\to \RR_{\ge 0}^P \times [-1,1]^Q 
\times \ti{\CC}_+^R\injto \RR^p\times \RR^{m-p}$$
obtained as the product of the maps
$$\begin{array}{ll}
\ti{\CC}_+^P \to \RR_{\ge 0}^P\,,&(z_i) \mapsto (r_i)\,, \\
\ti{\CC}_+^Q \to [-1, 1]^Q\,,&(r_j, \tau_j)\mapsto (\tau_j)\,,\\
\ti{\CC}_+^R \to \ti{\CC}_+^R\,, &\text{the identity map\,.}
\end{array}
$$
Note that $q=q_1\scirc pr$.
We thus have $\pi q^{-1}(A)=\pi_1 q_1^{-1}(A)$; we will show below
the allowability of $\pi_1 q_1^{-1}(A)$.

For a subset $J$ of $\{1, \cdots, n\}$,  the set 
$\BH_J\cap (\CC^P\times \CC_+^Q\times \CC_+^R )$ 
 will be just written $\BH_J$.
Let
$$\BHZ_J:=\{(z_i)\in \CC^P\times \CC_+^Q\times \CC_+^{R}\,|\,
z_i=0\,(i\in J)\,, \quad z_i\neq 0\, (i\in \{1, \cdots, n\} -J\,)\}$$
be an open set of $\BH_J$.
For $J=\emptyset$, $\BHZ_\emptyset$ is the complement of 
the union of $\BH_i$ for $i=1, \cdots, n$.
Clearly,  $\CC^P\times \CC_+^Q\times \CC_+^{R}$
 is the disjoint union of $\BHZ_J$, with $J$ varying over subsets of $\{1, \cdots, n\}$, and
$\BH_I=\cup_{J\supset I} \BHZ_J$. 

For a subset $J\subset \{1, \cdots, n\}$, let $J'=J\cap P$ and $J''=J\cap (Q\cup R)$.
Over $\BHZ_J$, the induced map $\pi_1: \pi_1^{-1}(\BHZ_J)\to \BHZ_J$ is a product bundle with fiber $[-1, 1]^{J''}$. 
Indeed the map $\pi: \ti{\CC}_+^Q\times \ti{\CC}_+^R\to \CC_+^Q\times \CC_+^R$ is a product 
bundle over the subset
$\{(z_i)\in \CC_+^Q\times \CC_+^R\, |\, z_i=0 \, (i\in J''),  z_i\neq 0\, (i\not\in J'')\,\}$
with fiber $[-1, 1]^{J''}$. 

Since $A$ is compact and $m$-admissible, for a subset $J$ of $\{1, \cdots, n\}$, one has
$\dim (A\cap \BH_J\cap \sq^n)\le m-2|J|$.
We will show  
$$\dim q_1\pi_1^{-1}(A\cap \BHZ_J) <m-|J'|\eqno{(\admissallow.a)}$$
for $J$ non-empty. Also we have 
$\dim q_1\pi_1^{-1}(A\cap \BHZ_\emptyset)\le m$, which 
is obvious since $\pi_1$ is an isomorphism over $A\cap \BHZ_\emptyset$.

For each  non-empty $J$, we can write
$$A\cap \BHZ_J= A_J'\cup A_J''$$
where $A_J'$, $A_J''$ are closed semi-algebraic sets of $\BHZ_J$
with
$$\dim A_J'\le m- 2|J| \quad\text{and}\quad A_J''\subset 
\BHZ_J\cap D \,.$$
Indeed one can take $A_J'$ to be the closure of $A\cap \BHZ_J\cap \sq^n$
in $\BHZ_J$ (with respect to the Euclidean
topology), and
$A''_J=A\cap \BHZ_J\cap D$.

Since $\pi_1$ is proper, continuous and semi-algebraic, 
$\pi_1^{-1}(A)$ is compact and semi-algebraic.
For a non-empty
subset $J\subset \{1, \cdots, n\}$ we have
$$\dim \pi_1^{-1}(A'_J)= \dim (A'_J) +|J'' | 
\le m-2|J|+|J''| =m-|J|-|J'|\,$$
(the last identity holds by $|J'|+|J''|= |J|$).
Thus 
 $\dim q_1\pi_1^{-1}(A'_J) <m-|J'|$.

As for $A_J''$, for each $t\in \{1, \cdots, n\}$, consider the set 
$A''_J\cap \{z_t=1\}$. 
If $t\in J$, then $A''_J\cap \{z_t=1\}$ is empty since $\BHZ_J\cap \{z_t=1\}=\emptyset$. 
If $t\not\in J$, then 
$$q_1\pi_1^{-1}(A''_J\cap \{z_t=1\})\subset
\left\{
\begin{array}{cl}
H_{J'}\cap \{r_t=1\}&\text{if $t\in P$\,,}\\
H_{J'}\cap\{\tau_t=0\}&\text{if $t\in Q$\,,}\\
H_{J'}\cap\{r_t=1\}\cap \{\tau_t=0\}&\text{if $t\in R$}\,,
\end{array}
\right.
$$
so in either case
$$\dim q_1\pi_1^{-1}(A''_J\cap \{z_t=1\})< \dim H_{J'} =m-|J'|\,.$$
Since $A''_J$ is the union of $A''_J\cap \{z_t=1\}$ for $t$, one has
$\dim q_1\pi_1^{-1}(A''_J) <m-|J'|$ for $J$ non-empty.
The inequality (\admissallow.a) hence follows.

Now, for $I$ a subset of $P$, we evaluate the dimension of the set
$q_1\pi_1^{-1}(A)\cap H_I$. 
By $q_1^{-1}(H_I)=\pi_1^{-1}(\BH_I)=\cup_{J\supset I} \pi_1^{-1}(\BHZ_J)$, 
one has
$$q_1\pi_1^{-1}(A)\cap H_I =\bigcup_{J\supset I} q_1\pi_1^{-1}(A\cap \BHZ_J)\,.
\eqno{(\admissallow.b)}$$
If $I\neq\emptyset$, for each $J$ containing $I$ (namely $J'\supset I$) we have
$$\dim q_1\pi_1^{-1}(A\cap \BHZ_J) <m-|J'|\le m-|I|\,,$$
so it follows that $\dim q_1\pi_1^{-1}(A)\cap H_I<m-|I|$. 
This shows that the set $q_1\pi_1^{-1}(A)$ is allowable with respect to $H_i$, $i\in P$.
\end{proof}

From this proposition and the discussion in (\Spolar)
we derive the following theorem.
\bigskip

\noindent (\convadmiss) {\bf Theorem.}\quad{\it 
 Let $A$ be a compact semi-algebraic set 
of $\CC^n$ which is $m$-admissible, and  
let  $\omega$ be an $(n, m-n)$-form on $\CC^n$ with
logarithmic singularities. 
Then the integral 
$$\int_A \omega$$
absolutely converges.

The same holds for an $m$-admissible compact semi-algebraic set of 
$(\BP^1)^n$ and a logarithmic $(n, m-n)$-form.
}\bigskip 

We turn to another related problem, and we show  
Theorem (\convadmisslink), which will be needed in the proof of the Cauchy formula.
\bigskip 

\noindent (\admissallowext) {\bf Proposition.}\quad{\it 
Let $n\ge 2$ and $A$ be  a compact semi-algebraic set 
of $\CC^n$ which is $m$-admissible. 
Let $\pi$ be a map as (\Spolar.b). 
Let $P, Q, R$ be a partition of $\{1, \cdots, n\}$ with $|P|+|Q|+2|R|=m$ and 
with $1\in R$, and $q: \ti{\CC}_+^n \to \RR^m$ be the associated map. 
Assume that 
$$A\subset \{|z_1|\ge |z_2|\}\,.$$
Then the set $q\pi^{-1}(A)$ is allowable with respect to 
$H_i=\{r_i=0\}$ for $i\in P\cup\{1\}$. 
}\smallskip 

{\bf Remark.}\quad   
(1) Compared with (\admissallow), the assumption is stronger and 
the conclusion is also stronger.  

(2) Without the assumption $A\subset \{|z_1|\ge |z_2|\}$,
the set $q \pi^{-1}(A)$ may not be allowable.

Indeed let $N=\{|z_1|\le 1\}\subset \CC$, 
let $B$ be a compact semi-algebraic set of dimension one of 
$\CC$ with $0\not\in B$, and define
$$A=N\times B\,;$$
then $A$ is compact 3-admissible semi-algebraic set of in $\CC^2$.
One easily verifies that the set $\ q \pi^{-1}(A)$ is not allowable
as its intersection with $\{r_1=0\}$ is of dimension 2.
\smallskip 

\begin{proof}  
We adapt the proof of Proposition (\admissallow), keeping the notation there.
With the maps $\pi_1$ and $q_1$ defined there, 
we will show the allowability of $q_1 \pi_1^{-1}(A)$. 
By the assumption on $A$, we have $A\cap \BH_1=A\cap \BH_{12}$, so 
$A\cap \BHZ_J=\empty$ if $1\in J$ but $2\not\in J$. 

We claim that for a non-empty subset $J$ of $\{1, \cdots, n\}$, the following
inequality holds:
$$\dim q_1\pi_1^{-1}(A\cap \BHZ_J) <m-|J\cap (P\cup \{1\})|\,.\eqno{(\admissallowext.a)}$$
We will derive this from the same inequalities with $A\cap \BHZ_J$ replaced by
$A'_J$ and $A''_J$.

We first show 
$$\dim \pi_1^{-1}(A'_J) <m-|J\cap (P\cup \{1\})|\,.\eqno{(\admissallowext.b)}$$
Indeed for $A'_J$ we showed in the proof of (\admissallow) the inequality
$$\dim \pi_1^{-1}(A'_J) \le m- |J|-|J'|\,.$$
If $1\not\in J$, then $J\cap (P\cup \{1\})=J'$ so (\admissallowext.c) holds. 
If $1\in J$ and $2\not\in J$, then $A\cap \BHZ_J$ is empty and the assertion
is trivial.  In case $\{1, 2\}\subset J$, then $|J\cap (P\cup \{1\})|=|J'|+1$ and 
$|J|>1$, hence
$$ m- |J|-|J'| < m-| J\cap (P\cup \{1\})|\,.$$

Next we show 
$$\dim q_1\pi_1^{-1}(A''_J) <m-|J\cap (P\cup \{1\})|\,.\eqno{(\admissallowext.c)}$$
If $t\in J$, then $A''_J\cap \{z_t=1\}$ is empty.
Assume thus $t\not\in J$.
If $1\not\in J$, this was shown in the proof of (\admissallow).
If $1\in J$ and $2\not\in J$, then $A\cap \BHZ_J=\emptyset$.
If $\{1, 2\}\subset J$ and $t\not\in J$,
then
$$q_1\pi_1^{-1}(A''_J\cap \{z_t=1\})\subset
\left\{
\begin{array}{cl}
H_{J'}\cap \{r_t=1\}\cap \{r_1=0\}&\text{if $t\in P$\,,}\\
H_{J'}\cap\{\tau_t=0\}\cap \{r_1=0\}&\text{if $t\in Q$\,,}\\
H_{J'}\cap\{r_t=1\}\cap \{\tau_t=0\}\cap \{r_1=0\}&\text{if $t\in R$}\,.
\end{array}
\right.
$$
Hence have $\dim  q_1\pi_1^{-1}(A''_J\cap \{z_t=1\})< \dim H_{J'}-1 =m-|J'|-1
=m-| J\cap (P\cup \{1\})|$, so (\admissallowext.c) holds.

For a non-empty subset $I$ of $P\cup\{1\}$, let $H_I=\cap_{i\in I} H_i$, and
 we evaluate the dimension  of
$q_1\pi_1^{-1}(A)\cap H_I$.
For each $J$ containing $I$, thus $J\cap (P\cup\{1\})\supset I$, it follows from 
(\admissallow.b) and (\admissallowext.a) that
$\dim q_1\pi_1^{-1}(A\cap \BHZ_J) <m-|I|$. 
Consequently $\dim q\pi_1^{-1}(A)\cap H_I<m-|I|$, concluding the proof.
\end{proof}

We state an analogous result for $A$ meeting $\BH_i$ only in $D$.
The proof is similar but simpler.
\bigskip 

\noindent (\admissallowextD) {\bf Proposition.}\quad{\it 
Let $A$ be  a compact semi-algebraic set 
of $\CC^n$ of dimension $\le m$ such that
$A\cap (\cup \BH_i)\subset D$.  
Let $\pi$ be a map as (\Spolar.b). 
Let $P, Q, R$ be a partition of $\{1, \cdots, n\}$ with $|P|+|Q|+2|R|=m$, 
and $q$ the associated map.
Then $q\pi^{-1}(A)$ is allowable with respect to 
$H_i=\{r_i=0\}$ for $i\in P\cup R$. 
}\smallskip 

\begin{proof}
We may assume $A\cap (\cup \BH_i)\subset \{z_t=1\}$ for some $t$.
Indeed, for a permutation $\sigma$ of $\{1, \ldots, n\}$, let 
$$R_\sigma=\{\,|z_{\sigma(1)}-1|\le |z_{\sigma(2)}-1|\le\cdots |z_{\sigma(n)}-1|\,\}\,.$$
Then $\CC^n$ is the union of $R_\sigma$ for permutations $\sigma$, and 
$D\cap R_\sigma\subset \{z_{\sigma(1)}=1\}$. 
The given $A$ is the union of $A\cap R_\sigma$, and one has
$(A\cap R_\sigma)\cap (\cup \BH_i)\subset \{z_{\sigma(1)}=1\}$.

For a non-empty subset $I$ of $P\cup R$, we evaluate the dimension 
of $q\pi^{-1}(A)\cap H_I$. 
If $t\in I$, the set is empty; thus we will assume $t\in I$. 
Then 
$$q\pi^{-1}(A)\cap H_I\subset
\left\{
\begin{array}{cl}
H_{I}\cap \{r_t=1\}&\text{if $t\in P$\,,}\\
H_{I}\cap\{\tau_t=0\}&\text{if $t\in Q$\,,}\\
H_{I}\cap\{r_t=1\}\cap \{\tau_t=0\}&\text{if $t\in R$}\,,
\end{array}
\right.
$$
and it follows that $\dim q\pi^{-1}(A)\cap H_I <\dim H_I$.
\end{proof}

\noindent (\convadmisslink) {\bf Theorem.}\quad{\it 
Suppose $n\ge 2$ and $A$ is a compact
semi-algebraic set 
of $\CC^n$ which is  $m$-admissible with respect to $\{\BH_i\}$.
For small $t>0$, let 
$$A_t =A\cap \{|z_1|=t \ge |z_2|\}\,, $$
which is a compact semi-algebraic set of dimension $\le m-1$.
Let $\omega$ be an $(n, m-n-1)$-form on $\CC^n$ with 
logarithmic singularities.
Then the integral 
$$\int_{A_t} |\omega|$$
converges to zero as $t\to 0$. 

The same holds for an $m$-admissible compact semi-algebraic set $A$ of 
$(\BP^1)^n$ and a logarithmic $(n, m-n)$-form $\omega$, taking 
$$A_t=A\cap \{|z_1|^{\eps_1}=t \ge |z_2|^{\eps_2}\}\,, $$
where $\eps_i\in \{-1, 1\}$ for $i=1, 2$.
}\smallskip 

\begin{proof} The case $A\subset (\BP^1)^n$ can be reduced to the affine case 
by an argument as in (\Spolar). We will thus assume $A\subset \CC^n$. 

One may replace $A$ with $A\cap \{|z_1|\ge |z_2|\}$, since the set
$A_t$ does not change.
So we will assume that $A$ is contained in the region $\{|z_1|\ge |z_2|\}$
and $A_t=A\cap \{|z_1|=t\}$. 

For $t>0$ small,  $A_t$ is a semi-algebraic set of dimension $\le m-1$
and $\dim(A_t\cap \cup \BH_i)\le m-2$,
so the 
integral $vol(A_t; \omega)=\int_{A_t}|\omega|$ is defined.
We will show that there exist $C, \al>0$ such that for $t>0$ small,
$$vol(A_t; \omega)\le Ct^\al\,.$$

Let $A^0=A-(\cup \BH_i)$  and $A_t^0=A_t\cap A^0$. 
Arguing as in (\Spolar), it is enough to show that, for each partition 
$(P, Q, R)$ of $\{1, \cdots, n\}$ with $|P|+|Q|+2|R|=m-1$, we have
an estimate of the form
$$vol(\pi^{-1}(A^0_t); \vphi_{P, Q, R})\le Ct^\al\,.\eqno{(\convadmisslink.a)}$$
If $1\in P\cup R$, then the form $\vphi_{P, Q, R}$ involves $dr_1$, so it restricts to zero on 
$A_t$ and the integral is zero.  Thus we may assume $1\in Q$.

Let $\bar{Q}= Q-\{1\}$ and $\bar{R}=R\cup\{1\}$; then $(P, \bar{Q}, \bar{R})$
is a partition of $\{1,\cdots, n\}$ with $|P|+|\bar{Q}|+2|\bar{R}|=m$.
We have the map
$$q: \ti{\CC}_+^n \to \RR^P\times \RR^{\bar Q}\times (\RR^2)^{\bar R}=\RR^m
\,$$
associated with $(P, \bar{Q}, \bar{R})$.
If $\vphi'_{P, Q, R}$ is the $(m-1)$-form on $\RR^m$
defined as (\Spolar.c),  $\vphi_{P, Q,R}=q^*\vphi'_{P, Q, R}$.

By (\admissallowext), $q\pi^{-1}(A)$ is allowable with respect to 
$H_i$, $i\in P\cup\{1\}$.  Thus by Corollary (\Thmnormallink.3) we have, with 
some $C, \al>0$,
$$vol(q\pi^{-1}(A)\cap \{r_1=t\};  \vphi'_{P, Q, R})\le C t^\al\,.\eqno{(\convadmisslink.b)}$$
Note that $q\pi^{-1}(A)\cap \{r_1=t\}=q\pi^{-1}(A_t)$.

Let $q: \pi^{-1}(A) \to \RR^m$ be the restriction of $q$; then we have a 
fiber square
$$\begin{array}{ccc}
\pi^{-1}(A_t)    &\hooklongrightarrow  &\pi^{-1}(A)  \\
\mapd{q} &     &\mapdr{q} \\
 \RR^{m-1}  &\overset{r_1=t}{\hooklongrightarrow} &\phantom{\,.}\RR^m\,.
 \end{array}
$$  
It follows that 
$\delta(q: \pi^{-1}(A_t) \to \RR^{m-1}) $ is bounded by 
$\delta:=\delta({q}: \pi^{-1}(A)\to \RR^m)$.

Applying \APropboundintegral to $q: \pi^{-1}(A_t)\to \RR^{m-1}$, we thus obtain
$$vol(\pi^{-1}(A_t); \vphi_{P, Q, R})\le \delta \cdot vol(q\pi^{-1}(A_t); \vphi'_{P, Q, R})\,.
\eqno{(\convadmisslink.c)}$$
The assertion (\convadmisslink.a) follows from (\convadmisslink.b) and (\convadmisslink.c).
\end{proof}

The proof of the following theorem is parallel, using
(\admissallowextD).
\bigskip 

\noindent (\convadmisslinkD) {\bf Theorem.}\quad{\it 
Let $A$ be as in 
(\admissallowextD), and 
$A_t:=A\cap \{|z_1|=t\}$. 
Then
for $\omega$ an $(n, m-n-1)$-form on $\CC^n$ with 
logarithmic singularities, the integral
$\int_{A_t} |\omega|$
converges to zero as $t\to 0$. 

The same holds for a compact semi-algebraic set $A$ of 
$(\BP^1)^n$ of dimension $\le m$ such that $A\cap (\cup \BH_I)\subset D$,
and a logarithmic $(n, m-n)$-form $\omega$, taking 
$$A_t=A\cap \{|z_1|^{\eps}=t\}\,, $$
where $\eps\in \{-1, 1\}$.
}\bigskip

{\bf References}
\smallskip

[Bl] Bloch, S.: The moving lemma for higher Chow groups, J. Algebraic
Geom.  3 (1994), 537--568.

[BCR] Bochnak, J., Coste, M., and Roy, M.-F. : Real algebraic geometry, Ergebnisse der Mathematik und ihrer Grenzgebiete vol. 36, Springer-Verlag, Berlin, 1998

[He] Herrera, M. E.: Integration on a semianalytic set. Bull. Soc. Math. France 94 (1966), 141-180.

[OS] Ohmoto, T. and Shiota, M. :
$C^1$-triangulations of semi-algebraic sets, preprint.
\bigskip 

\end{document}